\documentclass[aos,MSNbibl,nameyear,seceqn,dvips]{arximspdf}
\usepackage{graphicx}
\doi{10.1214/13-AOS1192} 
\volume{42}
\issue{1}
\pubyear{2014}
\firstpage{352}
\lastpage{381}

\makeatletter
\newcommand{\rrVert}{\Vert}
\newcommand{\rrvert}{\vert}
\newcommand{\llVert}{\Vert}
\newcommand{\llvert}{\vert}
\newtheorem{teo}{Theorem}[section]
\newtheorem{lem}[teo]{Lemma}
\newtheorem{prop}[teo]{Proposition}
\newproclaim{rem}[teo]{Remark}
\newproclaim{exa}{Example}
\def\T{\mathrm{T}}
\newcommand{\bfal}{\bolds{\alpha}}
\newcommand{\bfaa}{\mathbf{a}}
\newcommand{\bfaca}{\mathbf{A}}
\newcommand{\bfala}{\bolds{\alpha}}
\newcommand{\bfu}{\mathbf{u}}
\makeatother

\begin{document}
\begin{frontmatter}

\title{Anisotropic function estimation using multi-bandwidth Gaussian processes\thanksref{T1}}
\runtitle{Anisotropic function estimation}

\begin{aug}
\author[A]{\fnms{Anirban} \snm{Bhattacharya}\corref{}\ead[label=e1]{anirbanb@stat.tamu.edu}},
\author[B]{\fnms{Debdeep} \snm{Pati}\ead[label=e2]{debdeep@stat.fsu.edu}}
\and
\author[C]{\fnms{David} \snm{Dunson}\ead[label=e3]{dunson@stat.duke.edu}}
\runauthor{A. Bhattacharya, D. Pati and D. Dunson}
\affiliation{Texas A\&M University, Florida State University and Duke
University}
\address[A]{A. Bhattacharya\\
Department of Statistics\\
Texas A\&M University\\
College Station, Texas 77843\\
USA\\
\printead{e1}} 
\address[B]{D. Pati\\
Department of Statistics\\
Florida State University\\
Tallahassee, Florida 32306-4330\\
USA\\
\printead{e2}}
\address[C]{D. Dunson\\
Department of Statistical Science\\
Duke University\\
Durham, North Carolina 27708-0251\\
USA\\
\printead{e3}}
\end{aug}
\thankstext{T1}{Supported by Award Number R01ES017436 from the National
Institute of Environmental Health
Sciences and was conducted under the DARPA MSEE program.}


\received{\smonth{11} \syear{2011}}
\revised{\smonth{10} \syear{2013}}

%
\begin{abstract}
In nonparametric regression problems involving multiple
predictors, there is typically interest in estimating an anisotropic
multivariate regression
surface in the important predictors while discarding the unimportant ones.
Our focus is on defining a Bayesian procedure that leads to the minimax
optimal rate of posterior contraction (up to a log factor) adapting to
the unknown
dimension and anisotropic smoothness of the true surface. We propose
such an
approach based on a Gaussian process prior with dimension-specific
scalings, which
are assigned carefully-chosen hyperpriors.
We additionally show that using a homogenous Gaussian
process with a single bandwidth leads to a sub-optimal rate in
anisotropic cases.
\end{abstract}

%
\begin{keyword}[class=AMS]
\kwd[Primary ]{62G07}
\kwd{62G20}
\kwd[; secondary ]{60K35}
\end{keyword}
\begin{keyword}
\kwd{Adaptive}
\kwd{anisotropic}
\kwd{Bayesian nonparametrics}
\kwd{function estimation}
\kwd{Gaussian process}
\kwd{rate of convergence}
\end{keyword}

\end{frontmatter}

\setcounter{footnote}{1}
\section{Introduction}\label{secintro}
Gaussian processes [\citet{rasmussen2004gaussian}, \citet
{van2008reproducing}] are widely used as priors on functions due to
tractable posterior computation and attractive theoretical properties.
The law of a mean zero Gaussian process (GP) $W_t$ is entirely
characterized by its covariance kernel $c(s,t) = E(W_s W_t)$. A squared
exponential covariance kernel given by $c(s, t) = \exp(-a \llVert
s-t\rrVert ^2)$ is commonly used in the literature.

Given $n$ independent observations, the
optimal rate of estimation of a $d$-variable function that is only
known to be $\alpha$-smooth is $n^{-\alpha/(2 \alpha+ d)}$
[\citet{stone1982optimal}]. The quality of estimation thus
improves with increasing smoothness of the ``true'' function while it
deteriorates with increase in dimensionality. In practice, the
smoothness $\alpha$ is typically unknown and one would like an
estimation procedure that adapts to any possible $\alpha>0$.
Accordingly, a lot of effort has been employed to develop adaptive
estimation methods that are rate-optimal for every regularity level of
the unknown function.

The literature on adaptive estimation in a minimax setting was
initiated by Lepski{\u\i} in a series of papers [\citeauthor{lepski1990problem} (\citeyear{lepski1990problem,lepski1991asymptotic,lepski1992asymptotic})];
see also \citet{birge2001alternative}
for a discussion on this topic. We also refer the reader to \citet
{hoffmann2002random}, which contains an extensive list of developments
in the frequentist literature on adaptive estimation. There is a
growing literature on Bayesian adaptation over the last decade.
Previous works include \citet{belitser2003adaptive},
\citeauthor{ghosal2003bayesian} (\citeyear{ghosal2003bayesian,ghosal2008nonparametric}),
\citet{huang2004convergence}, \citet
{scricciolo2006convergence}, \citet{rousseau2010rates},
\citet{kruijer2010adaptive},
\citet{de2010adaptive}, \citet{shen2011adaptive}.

A key idea in frequentist adaptive estimation is to narrow down the
search for an ``optimal'' estimator within a class of estimators
indexed by a smoothness or bandwidth parameter, and make a data-driven
choice to select the proper bandwidth. In a Bayesian context, one would
place a prior on the bandwidth parameter and model-average across
different values of the bandwidth through the posterior distribution.
The parameter $a$ in the squared-exponential covariance kernel $c$
plays the role of a scaling or inverse bandwidth. \citet
{van2009adaptive} showed that with a gamma prior on $a^d$, one obtains
the minimax rate of posterior contraction $n^{-\alpha/(2 \alpha+ d)}$
up to a logarithmic factor for $\alpha$-smooth functions adaptively
over all $\alpha> 0$.

In most multivariate applications, the isotropic smoothness assumption
seems too restrictive. Potentially, one can incorporate a separate
scaling variable $a_j$ for each dimension using the covariance kernel
$c(s,t) =\break  \exp(-\sum_{j=1}^d a_j |s_j-t_j|^2)$, intuitively enabling
better approximation of\break  anisotropic functions. Such kernels, going by
the name \emph{automatic relevance determination} (ARD), have been
heavily used in the machine learning community; see, for example, \citet{rasmussen2004gaussian} and references therein.
\citet{zou2010nonparametric} and \citet{savitsky2011variable} recently
considered such a model, with point mass mixture priors on the $a_j$'s.
Although this is an attractive approach with encouraging empirical
performance, there has not been any theoretical studies of asymptotic
properties of related models in a Bayesian framework.


In the frequentist literature, minimax rates of convergence in
anisotropic Sobolev, Besov and H\"{o}lder spaces have been studied in
\citet{ibragimov1981statistical}, \citet{nussbaum1985spline}, \citet{birge1986estimating}, with
adaptive estimation procedures developed in \citet{barron1999risk},
\citet{kerkyacharian2001nonlinear}, \citet{hoffmann2002random}, \citet{klutchnikoff2005adaptive}
among others. The traditional way of dealing with anisotropy is to
employ a separate bandwidth or scaling parameter for the different
dimensions, and choose an optimal combination of scales in a
data-driven way. However, the multidimensional nature of the problem
makes the optimal bandwidth selection difficult compared to the
isotropic case, as there is no natural ordering among the estimators
with multiple bandwidths [\citet{lepski1999adaptive}].

It is known [\citet{hoffmann2002random}] that the minimax rate of
convergence for a function with smoothness $\alpha_i$ along the $i$th
dimension is given by $n^{-\alpha_0/(2\alpha_0 + 1)}$, where $\alpha
_0^{-1} = \sum_{i=1}^d \alpha_i^{-1}$ is an \emph{exponent of global
smoothness} [\citet{birge1986estimating}]. When $\alpha_i =
\alpha$ for all $i =1, \ldots, d$, one reduces back to the optimal
rate for isotropic classes. On the contrary, if the true function
belongs to an anisotropic class, the assumption of isotropy would lead
to loss of efficiency which would be more and more accentuated in
higher dimensions.
In addition, if the true function depends on a subset of coordinates $I
= \{i_1, \ldots, i_{d_0}\} \subset\{1, \ldots, d\}$ for some $1 \leq
d_0 \leq d$, the minimax rate would further improve to $n^{-\alpha
_{0I}/(2 \alpha_{0I} + 1)}$, with $\alpha_{0I}^{-1} = \sum_{j \in I}
\alpha_j^{-1}$.

The objective of this article is to study whether one can fully adapt
to this larger class of functions in a Bayesian framework using
dimension-specific rescalings of a homogenous Gaussian process,
referred to as a multi-bandwidth Gaussian process from now on. We
answer the question in the affirmative to establish rate adaptiveness
of the posterior distribution in a variety of settings involving a
multi-bandwidth Gaussian process through a novel prior specification on
the vector of bandwidths. For simplicity of exposition, we initially
study the problem in two parts: (i) adaptive estimation over
anisotropic H\"{o}lder functions of $d$ arguments, and (ii) adaptive
estimation over functions that can possibly depend on fewer coordinates
and have isotropic H\"{o}lder smoothness over the remaining coordinates.
The proposed prior specification for the two cases above are
intuitively interpretable and can be easily connected to prescribe a
unified prior leading to adaptivity over (i) and (ii) combined.

Although our prior specification involving dimension-specific bandwidth
parameters leads to adaptivity, a stronger result is required to
conclude that a single bandwidth would be inadequate for the above
classes of functions. We prove that the optimal prior choice in the
isotropic case leads to a sub-optimal convergence rate if the true
function has anisotropic smoothness by obtaining a lower bound on the
posterior contraction rate. Previous results on posterior lower bounds
in nonparametric problems include \citet{castillo2008lower}, \citet{van2011information}.

The remaining paper is organized as follows. In Section~\ref{secnotn}, we introduce relevant notations and conventions used
throughout the paper. The multi-bandwidth Gaussian process is
introduced in Section~\ref{secmainres}. Sections~\ref{secaniso}
and~\ref{secadap} discuss the main developments with applications to
anisotropic Gaussian process mean regression and logistic Gaussian
process density estimation described in Section~\ref{subsecspecific}.
Section~\ref{seclowerbounds} establishes the necessity of the
multi-bandwidth Gaussian process by showing a~lower-bound result.
In Sections~\ref{secproprescaled} and~\ref{ssecauxlb}, we study
various properties of rescaled Gaussian processes which are crucially
used in the proofs of the main theorems in
Section~\ref{secmainproofs}.\looseness=1

\section{Preliminaries}\label{secnotn}
To keep the notation clean, we shall only use boldface for $\mathbf
{a}, \mathbf{b}$ and $\bfal$ to denote vectors.
We shall make frequent use of the following multi-index notations. For
vectors $\mathbf{a}, \mathbf{b}$ $\in\mathbb{R}^d$, let $\mathbf
{a}. = \sum_{j=1}^d a_j, \mathbf{a}^* = \prod_{j=1}^d a_j, \mathbf
{a}! = \prod_{j=1}^d a_j !, \bar{\mathbf{a}} = \max_j{a_j},
\underline{\mathbf{a}} = \min_j{a_j}$,\vspace*{-1pt} $\mathbf{a}./\mathbf{b}=
(a_1/b_1, \ldots, a_d/b_d)^{\T}$, $\mathbf{a}\cdot\mathbf{b}=
(a_1b_1, \ldots, a_db_d)^{\T}$, $\mathbf{a}^{\mathbf{b}} = \prod_{j=1}^d a_j^{b_j}$. Denote\vspace*{2pt} $\mathbf{a}\leq\mathbf{b}$ if $a_j \leq
b_j$ for all $j = 1, \ldots, d$. For $n = (n_1, \ldots, n_d)$, let
$D^{n} f$ denote the mixed partial derivatives of order $(n_1, \ldots,
n_d)$ of $f$.

Let $C[0, 1]^d$ and $C^{\beta}[0, 1]^d$ denote the space of all
continuous functions and the H\"{o}lder space of $\beta$-smooth
functions $f\dvtx  [0, 1]^d \to\mathbb{R}$, respectively, endowed with the
supremum norm $\llVert  f\rrVert _{\infty} = \sup_{ t \in[0,
1]^d } \llvert  f(t) \rrvert $. For $\beta> 0$, the H\"{o}lder
space $C^{\beta}[0, 1]^d$ consists of functions $f \in C[0, 1]^d$ that
have bounded mixed partial derivatives up to order $\lfloor\beta
\rfloor$, with the partial derivatives of order $\lfloor\beta\rfloor
$ being Lipschitz continuous of order $\beta- \lfloor\beta\rfloor$.
Also, denote by $H^{\beta}[0, 1]^d$ the Sobolev space of functions $f\dvtx
[0, 1]^d \to\mathbb{R}$ that are restrictions of a function $f\dvtx
\mathbb{R}^d \to\mathbb{R}$ with Fourier transform $\hat{f}(\lambda
) = (2\pi)^{-d}\int e^{i(\lambda, t)} f(t) \,dt$ such that
\[
\int\bigl(1+ \llVert \lambda\rrVert ^2\bigr)^{\beta}\bigl
\llvert \hat {f}(\lambda)\bigr\rrvert ^2 \,d\lambda< \infty.
\]
Note that using the above convention, the inverse Fourier transform
$f(t) = \int e^{-i(\lambda, t)} \hat{f}(\lambda) \,d\lambda$.
Next, we define an anisotropic H\"{o}lder class of functions previously
used in \citet{barron1999risk}, \citet{klutchnikoff2005adaptive}. For a function
$f \in C[0, 1]^d$, $x \in[0, 1]^d$, and $1 \leq i \leq d$, let $f_i(
\cdot\mid x)$ denote the univariate function $y \mapsto f(x_1, \ldots, x_{i-1}, y, x_{i+1}, \ldots, x_d)$. For a vector of positive numbers
$\bfal= (\alpha_1, \ldots, \alpha_d)$, the anisotropic H\"{o}lder
space $C^{\bfala}[0, 1]^d$ consists of functions $f$ which satisfy,
for some $L > 0$,
%
%
\begin{equation}
\label{eqanisocond1} \max_{1 \leq i \leq n} \sup_{x \in[0, 1]^d} \sum
_{j=0}^{\lfloor
\alpha_i \rfloor} \bigl\llVert D^j
f_i(\cdot\mid x)\bigr\rrVert _{\infty} \leq L
\end{equation}
and, for any $y \in[0,1]$, $h$ small such that $y+h \in[0,1]$ and for
all $1\leq i \leq d$,
%
%
\begin{equation}
\label{eqanisocond2} \sup_{x \in[0,1]^d} \bigl\llVert D^{\lfloor\alpha_i \rfloor}
f_i(y+h \mid x) - D^{\lfloor\alpha_i \rfloor} f_i(y \mid x)\bigr
\rrVert _{\infty} \leq L\llvert h\rrvert ^{\alpha_i - \lfloor
\alpha_i \rfloor}.
\end{equation}

For $t \in\mathbb{R}^d$ and a subset $I \subset\{1, \ldots, d\}$ of
size $\llvert  I\rrvert = \tilde{d}$ with $1 \leq\tilde{d}
\leq d$, let $t_I$~denote the vector of size $\tilde{d}$ consisting of
the coordinates $(t_j\dvtx  j \in I)$. Let $C[0, 1]^I$ denote the subset of
$C[0, 1]^d$ consisting of functions $f$ such that $f(t) = g(t_I)$ for
some function $g \in C[0, 1]^{\tilde{d}}$.
Also, let $C^{\bfala}[0, 1]^I$ denote the subset of $C^{\bfala}[0,
1]^d$ consisting of functions $f$ such that $f(t) = g(t_I)$ for some
function $g \in C^{\bfal_I}[0, 1]^{\tilde{d}}$.

The $\varepsilon$-covering number $N(\varepsilon,S,d)$ of a semimetric
space $S$ relative to the semimetric $d$ is the minimal number of balls
of radius $\varepsilon$ needed to cover $S$. The logarithm of the
covering number is referred to as the entropy.

We write ``$\precsim$'' for inequality up to a constant multiple. Let
$\phi(x) = (2\pi)^{-1/2}\*\exp(-x^2/2)$ denote the standard
normal density, and let $\phi_{\sigma}(x) = (1/\sigma) \phi
(x/\sigma)$. Let an asterisk denote a convolution, for example, $(\phi
_{\sigma} * f)(y) = \int\phi_{\sigma}(y - x)f(x)\,dx$. 

Let $\mathbb{R}_+$ denote the set of nonnegative real numbers and let
$\mathbb{R}_*$ denote positive reals. Denote by $\mathcal{S}_{d-1}$
the $(d-1)$-dimensional simplex $\{ x \in\mathbb{R}^d\dvtx  x_i \geq0,
1\leq i \leq d, \sum_{i=1}^d x_i = 1\}$.

Unless otherwise stated, $c, C, C', C_1, C_2, \ldots, K, K_1, K_2,
\ldots$ shall denote global constants irrelevant to our purpose.
\section{Main results}\label{secmainres}
Let $W = \{W_t\dvtx  t \in[0,1]^d\}$ be a centered homogeneous Gaussian
process with covariance function $\mbox{E}(W_s W_t) = c(s-t)$. A
detailed review of the facts on Gaussian processes relevant to the
present application can be found in \citet{van2008reproducing}. If $c\dvtx
\mathbb{R}^d \to\mathbb{R}$ is continuous, by Bochner's theorem,
there exists a finite positive measure $\nu$ on $\mathbb{R}^d$,
called the spectral measure of $W$, such that
\[
c(t) = \int_{\mathbb{R}^d} e^{- i (\lambda, t)} \nu(d\lambda),
\]
where for $u, v \in\mathbb{C}^d$, $(u, v)$ denotes the complex inner product.
As in \citet{van2009adaptive}, we shall restrict ourselves to
processes with spectral measure $\nu$ having subexponential tails,
that is, for some $\delta> 0$,
%
%
\begin{equation}
\label{eqsubexp} \int e^{ \delta\llVert \lambda\rrVert } \nu(d\lambda) < \infty.
\end{equation}
The spectral measure $\nu$ of a squared exponential covariance kernel
with $c(t) = \exp(- \llVert  t\rrVert ^2)$ has a\vspace*{1pt} density
w.r.t. the Lebesgue measure given by $f(\lambda) = 1/(2^d \pi
^{d/2})\exp(- \llVert \lambda\rrVert ^2/4)$ which clearly
satisfies (\ref{eqsubexp}).

Let $\mathbb{H}$ be the RKHS of $W$; see \citet{van2008reproducing}
for a~review of relevant facts. \citet{van2008rates} showed that the
rate of posterior contraction with a Gaussian process prior $W$, using
a metric where appropriate testing is possible, is determined by the
behavior of the concentration function $\phi_{w_0}(\varepsilon)$ for
$\varepsilon$ close to zero, where
%
%
\begin{equation}
\phi_{w_0}(\varepsilon) = \inf_{h \in\mathbb{H}\dvtx  \llVert  h -
w_0\rrVert _{\infty} \leq\varepsilon} \frac{\llVert  h\rrVert ^2_{\mathbb{H}}}{2}
- \log P  \bigl(\llVert W\rrVert _{\infty
} \leq\varepsilon\bigr).
\end{equation}
We tacitly assume that there is a given statistical problem where the
true parameter $g_0$ is a known function of $w_0$.
\citet{van2009adaptive} studied rescaled Gaussian processes $W^A
= \{W_{At}\dvtx  t \in[0,1]^d\}$ for a real positive random variable $A$
stochastically independent of $W$, and showed that with a Gamma prior
on $A^d$, one obtains the minimax-optimal rate of convergence
$n^{-\alpha/(2\alpha+ d)}$ (up to a logarithmic factor) for $\alpha
$-smooth functions. Since their prior specification does not involve
the unknown smoothness $\alpha$, the procedure is fully adaptive.

The key result of \citet{van2009adaptive} was to construct sets $B_n
\subset C[0, 1]^d$ so that given $\alpha> 0$, a function $w_0 \in
C^{\alpha}[0, 1]^d$, and a constant \mbox{$C > 1$}, there exists a constant
$D > 0$ such that, for every sufficiently large $n$,
%
%
\begin{eqnarray}
\log N\bigl(\bar{\varepsilon}_n, B_n, \llVert \cdot
\rrVert _{\infty}\bigr) &\leq& D n \bar{\varepsilon}_n^2,
\label{eqsieveentropy}
\\
\mathrm{P}\bigl(W^{A} \notin B_n\bigr) &\leq&
e^{-Cn \varepsilon_n^2}, \label{eqsievecompl}
\\
\mathrm{P}\bigl(\bigl\llVert W^{A} - w_0\bigr\rrVert
_{\infty} \leq \varepsilon_n\bigr) &\geq& e^{- n \varepsilon_n^2}
\label{eqnoncentered}
\end{eqnarray}
with $\varepsilon_n = n^{-\alpha/(2 \alpha+ 1)} (\log n)^{\kappa_1},
\bar{\varepsilon}_n = n^{-\alpha/(2 \alpha+ 1)} (\log n)^{\kappa_2}$
for constants $\kappa_1, \kappa_2 > 0$.

In this article, we shall study multi-bandwidth Gaussian processes of
the form $\{ W_t^{\mathbf{a}} = W_{\mathbf{a}\cdot t}\dvtx    t \in
[0,1]^d\}$ for a vector of rescalings (or inverse-bandwidths) $\mathbf
{a}= (a_1, \ldots, a_d)^{\T}$ with $a_j > 0$ for all $j = 1, \ldots, d$.
We consider two function classes defined in Section~\ref{secnotn}:
\begin{longlist}[(ii)]
\item[(i)] H\"{o}lder class of functions $C^{\bfala}[0, 1]^d$ with
anisotropic smoothness (\mbox{$\bfal\!\in\!\mathbb{R}_*^d$}).
\item[(ii)] H\"{o}lder class of functions $C^{\alpha}[0, 1]^I$ with
isotropic smoothness that can possibly depend on fewer dimensions
($\alpha> 0$ and $I \subset\{1, \ldots, d\}$).
\end{longlist}
For a continuous function in the support of a Gaussian process, the
probability assigned to a sup-norm neighborhood of the function is
controlled by centered small ball probability and how well the function
can be approximated from the RKHS of
the process [Section~5 of \citet{van2008reproducing}]. With the target
class of functions as in (i)~or~(ii), a single scaling
seems inadequate and it is intuitively appealing to introduce multiple
bandwidth parameters to enlarge the RKHS and facilitate improved
approximation from the RKHS.

The main technical challenge for adaptation in our setting is to find a
joint prior on $\mathbf{a}$ and devise sets $B_n$ (henceforth called
\emph{sieves}) so that (\ref{eqsieveentropy})--(\ref{eqnoncentered}) are satisfied with $w_0$ in the above function
classes (i)--(ii) and $\varepsilon_n$ being the optimal rate
of convergence for the same. To that end, we propose a novel class of
joint priors on the rescaling vector $\mathbf{a}$ that leads to
adaptation over function classes (i)~and~(ii) in
Sections~\ref{secaniso} and~\ref{secadap}, respectively.
Connections between the two prior choices are discussed and a unified
framework is prescribed for the function class $ \{C^{\bfala}[0,
1]^{I}\dvtx  \bfal\in\mathbb{R}_*^d, I \subset\{1, \ldots, d\}  \}$
combining (i) and (ii).

The construction of the sieves $B_n$ are laid out in Section~\ref{secmainproofs}. With such $B_n$, one can use standard results to
establish adaptive minimax rate of convergence in various statistical
settings; refer to the discussion following Theorem 3.1 in \citet{van2009adaptive}.
Some such specific applications are described in Section~\ref{subsecspecific}.

\subsection{Adaptive estimation of anisotropic functions}\label{secaniso}
Let $\mathbf{A} = (A_1, \ldots, A_d)^{\T}$ be a random vector in
$\mathbb{R}^d$ with each $A_j$ a nonnegative random variable
stochastically independent of $W$. We can then define a scaled process
$W^{\bfaca} = \{W_{\mathbf{A}\cdot t}\dvtx  t \in[0,1]^d\}$, to be
interpreted as a Borel measurable map in $C[0, 1]^d$ equipped with the
sup-norm $\llVert \cdot\rrVert _{\infty}$. The basic idea
here is to scale the different dimensions by different amounts so that
the resulting process becomes suitable for approximating functions
having different smoothness along the different coordinate axes.

We shall define a joint distribution on $\mathbf{A}$ induced through
the following hierarchical specification. Let $\Theta= (\Theta_1,
\ldots, \Theta_d)$ denote a random vector with a density supported on
the simplex $\mathcal{S}_{d-1}$.
\begin{longlist}[(PA2)]
\item[(PA1)] Draw\vspace*{1pt} $\Theta\sim\operatorname{Dir}(\beta_1, \ldots, \beta
_d)$ for some $\beta= (\beta_1, \ldots, \beta_d)$.
\item[(PA2)] Given $\Theta= \theta$, draw $A_j^{1/\theta_j} \sim g$
independently, where $g$ is a density on the positive real line satisfying
%
%
\begin{equation}
\label{eqgdens3.6} C_1 x^p \exp\bigl(-D_1 x
\log^q x\bigr) \leq g(x) \leq C_2 x^p \exp
\bigl(-D_2 x \log^q x\bigr)
\end{equation}
for positive constants $C_1, C_2, D_1, D_2$ and nonnegative constants
$p$, $q$ and every sufficiently large $x > 0$.
\end{longlist}
In particular, $g$ corresponds to a $\operatorname{gamma}(b_1, b_2)$
distribution if $b_1 = p+1$, $D_1=D_2=b_2$, $q = 0$, $C_1=C_2$ in (\ref{eqgdens3.6}). For notational simplicity, we shall assume $g$ to be
$\operatorname{gamma}(b_1, b_2)$ from now on, noting that the main results
would all hold for the general form of $g$ above.

Let $\pi_{\bfaca}$ denote the joint prior on $\mathbf{A}$ induced
through (PA1)--(PA2), so that $\pi_{\bfaca}(\mathbf{a})
= \int\prod_{j=1}^d \pi(a_j \mid\theta_j) \,d\pi(\theta)$. We now
state our main theorem for the anisotropic smoothness class in~(i), with a detailed proof provided in Section~\ref{secmainproofs}.
%
%
\begin{teo}\label{teoanisomain}
Let $W$ be a centered homogeneous Gaussian random field on $\mathbb
{R}^d$ with spectral measure $\nu$ that satisfies (\ref{eqsubexp})
and let $W^{\bfaca}$ denote the multi-bandwidth process with $\mathbf
{A}\sim\pi_{\bfaca}$ as in \textup{(PA1)--(PA2)}.
Let $\bfal= (\alpha_1, \ldots, \alpha_d)$ be a vector of positive
numbers and $\alpha_0 = (\sum_{i = 1}^d \alpha_i^{-1})^{-1}$.
Suppose $w_0$ belongs to the anisotropic H\"{o}lder space $C^{\bfala
}[0, 1]^d$. Then for every constant $C > 1$, there exist Borel
measurable subsets $B_n$ of $C[0, 1]^d$ and a constant $D > 0$ such
that, for every sufficiently large $n$, the conditions (\ref
{eqsieveentropy})--(\ref{eqnoncentered}) are\vspace*{1pt} satisfied by
$W^{\bfaca}$ with $\varepsilon_n = n^{-\alpha_0/(2 \alpha_0 + 1)}
(\log n)^{\kappa_1}, \bar{\varepsilon}_n = n^{-\alpha_0/(2 \alpha_0 +
1)} (\log n)^{\kappa_2}$ for constants $\kappa_1, \kappa_2 > 0$.
\end{teo}

\subsection{Adaptive dimension reduction}\label{secadap}

We next consider the smoothness class in (ii), namely $C^{\alpha
}[0, 1]^I$ for $I \subset\{1, \ldots, d\}$ and $\alpha> 0$. If the
true function has isotropic smoothness on the dimensions it depends on,
it is intuitively clear that one does not need a separate scaling for
each of the dimensions. Indeed, had\vadjust{\goodbreak} we known the true coordinates $I
\subset\{1, \ldots, d\}$, we could have only scaled the dimensions in
$I$ by a positive random variable $A$, and a slight modification of the
results in \citet{van2009adaptive} would imply that a gamma prior on
$A^{\llvert  I\rrvert }$ would lead to adaptation.

With that motivation, consider a joint prior $\pi_{\bfaca}$ on
$\mathbf{A}$ induced through the following hierarchical scheme:
\begin{longlist}[(PD3)]
\item[(PD1)] draw $\tilde{d}$ uniformly on $\{0, \ldots, d\}$,
\item[(PD2)] given $\tilde{d}$, draw a subset $S \subset\{1, \ldots, d\}$ with $|S| = \tilde{d}$ uniformly from all subsets of size
$\tilde{d}$,
\item[(PD3)] generate a pair of random variables $(A, B)$ with
$A^{\tilde{d}} \sim\operatorname{gamma}(b_1, b_2)$ and $B$ drawn from a
fixed compactly supported distribution,
\item[(PD4)] set $A_j = A$ for $j \in S$ and $A_j = B$ for $j \notin S$.
\end{longlist}
In particular, one can fix $B$ to be any constant $c > 0$ in (PD3),
which corresponds to the $\delta_c(\cdot)$ distribution.
We next state our main result on adaptive dimension reduction. The
proof of the following Theorem~\ref{teodimredmain} has elements in
common with the proof of the Theorem~\ref{teoanisomain}, and hence
only a sketch of the proof is provided in Section~\ref{secmainproofs}. 

%
%
\begin{teo}\label{teodimredmain}
Let $W$ be a centered homogeneous Gaussian random field on $\mathbb
{R}^d$ with spectral measure $\nu$ that satisfies (\ref{eqsubexp})
and let $W^{\bfaca}$ denote the multi-bandwidth process with $\mathbf
{A}\sim\pi_{\bfaca}$ as in \textup{(PD1)--(PD4)}.
Suppose $w_0$ belongs to the H\"{o}lder space $C^{\alpha}[0, 1]^I$ for
some subset $I$ of $\{1, \ldots, d\}$ and $\alpha> 0$. Then for\vspace*{1pt} every
constant $C > 1$, there exist Borel measurable subsets $B_n$ of $C[0,
1]^d$ and a constant $D > 0$ such that, for every sufficiently large
$n$, the\vspace*{1pt} conditions (\ref{eqsieveentropy})--(\ref{eqnoncentered})
are satisfied by $W^{\bfaca}$ with $\varepsilon_n = n^{-\alpha/(2
\alpha+ d_0)} (\log n)^{\kappa_1}, \bar{\varepsilon}_n = n^{-\alpha
/(2 \alpha+ d_0)} (\log n)^{\kappa_2}$ for constants $\kappa_1,
\kappa_2 > 0$ and $d_0 = \llvert  I\rrvert $.
\end{teo}

%

\subsection{Connections between cases \textup{(i)} and \textup{(ii)}}\label{secconnec}

The joint distributions on $\mathbf{A}$ specified in (PA1)--(PA2) and
(PD1)--(PD4) are closely connected. To begin
with, note that if we set $A_j = A, \theta_j = 1/d, j = 1, \ldots, d$
in \textup{(PA1)--(PA2)}, one reduces to a gamma prior on $A^d$;
the optimal prior choice in the isotropic case [\citet
{van2009adaptive}]. In the anisotropic case, our proposed prior can be
motivated as follows. Recall that the purpose of rescaling is to
traverse the sample paths of an infinite smooth stochastic process on a
larger domain to make it more suitable for less smooth functions. If
the true function has anisotropic smoothness, then we would like to
stretch those directions more where the function is less smooth. For
smaller $\theta_j$'s, the marginal distribution of $a_j$ has lighter
tails compared to larger values of $\theta_j$. We would thus like
$\theta_j$ to assume smaller values for the directions $j$ where the
function is more smooth and larger values corresponding to the less
smooth directions. Without further constraints on $\theta$, it is not
possible to separate the scale of $\mathbf{A}$ from $\theta$. This
motivates us to constrain $\theta$ to the simplex which serves as a
weak identifiability condition.\vadjust{\goodbreak}

In the limit as $\theta_j \to0$, the distribution of $a_j$ converges
to a point mass at zero. Accordingly, if the true function does not
depend on a set of $(d-d^*)$ dimensions, we would set $\theta_j = 0$
for those dimensions and choose the remaining $\theta_j$'s from a
$d^*-1$-dimensional simplex. In particular, if the function has
isotropic smoothness in the remaining $d^*$ coordinates, one can simply
choose $\theta_j = 1/d^*$ for those dimensions. This reduces to our
prior choice in \textup{(PD1)--(PD4)}.

The dimensionality reduction in Section~\ref{secadap} deals with
finitely many models and can be alternatively studied as a model
selection problem [\citet{ghosal2008nonparametric}]. However, the
connection established above allows us to treat anisotropy and
dimension reduction under a single framework with the dimension
reduction paradigm recognized as a limiting case of the anisotropic
framework. We exploit this connection to propose a unified framework
for adaptively estimating functions which possibly depend on fewer
coordinates and have anisotropic smoothness in the remaining ones, that
is, functions in $C^{\bfala}[0,1]^I$ for $\bfal\in\mathbb{R}_*^d$
and $I \subset\{1, \ldots, d\}$. In particular, we continue to
consider rescaled Gaussian processes $W^{\bfaca}$ with the following
prior on $\mathbf{A}$:
\begin{longlist}[(P4)]
\item[(P1)] draw $\tilde{d}$ uniformly on $\{0, \ldots, d\}$,
\item[(P2)] given $\tilde{d}$, draw a subset $S \subset\{1, \ldots,
d\}$ with $|S| = \tilde{d}$ uniformly from all subsets of size $\tilde{d}$,
\item[(P3)] draw $\Theta= (\theta_1, \ldots, \theta_{\tilde{d}})
\in\mathcal{S}_{\tilde{d}-1}$ from a $\operatorname{Dir}(\beta_1, \ldots,
\beta_{\tilde{d}})$ distribution,
\item[(P4)] given $S$ and $\Theta= \theta$, draw $A_j^{1/\theta_j}
\sim\operatorname{gamma}(b_1, b_2)$ independently for $j \in S$, and fix $A_j
= c$ for some $c > 0$ and $j \notin S$.
\end{longlist}
A unification of Theorems~\ref{teoanisomain} and~\ref{teodimredmain} is provided in the following Theorem~\ref{teocombinedmain}, with the proof omitted as it is similar to the
previous theorems.
%
%
\begin{teo}\label{teocombinedmain}
Consider $W^{\bfaca}$ with $\mathbf{A}\sim\pi_{\bfaca}$ as in
\textup{(P1)--(P4)}.
Suppose $w_0$ belongs to $C^{\bfala}[0, 1]^I$ for some subset $I$ of
$\{1, \ldots, d\}$ and $\bfal\in\mathbb{R}_*^d$. Let $\alpha
_{0I}^{-1} = \sum_{j \in I} \alpha_j^{-1}$. Then for every constant
$C > 1$, the conclusions of Theorem~\ref{teoanisomain}~are satisfied
by $W^{\bfaca}$ with $\varepsilon_n = n^{-\alpha_{0I}/(2 \alpha_{0I} +
1)} (\log n)^{\kappa_1}, \bar{\varepsilon}_n =\break  n^{-\alpha_{0I}/(2
\alpha_{0I} + 1)} (\log n)^{\kappa_2}$ for constants $\kappa_1,
\kappa_2 > 0$.
\end{teo}
In Theorem~\ref{teocombinedmain}, the exponents of the logarithmic
terms increase linearly with the dimension and decrease with $\alpha
_{0I}$. In particular, $\kappa_1$ and $\kappa_2$ can be estimated by
$\frac{d+1}{\alpha_{0I} +2}$ and $\frac{(\alpha_{0I} +
3)(d+1)}{2(\alpha_{0I} + 2)}$, respectively.

\subsection{Rates of convergence in specific settings}\label{subsecspecific}
Theorem~\ref{teocombinedmain} is in the same spirit as Theorem 3.1
of \citet{van2009adaptive} [see also Theorem 2.2 of \citet{de2010adaptive}] and can be used to derive rates of posterior
contraction in a variety\vadjust{\goodbreak} of statistical problems involving Gaussian
random fields. We shall consider a couple of specific problems with the
message that similar results can be obtained for a large class of problems.

We first consider a regression problem where given independent response
variable $y_i$ and covariates $x_i \in[0,1]^d$, the response is
modeled as random perturbations around a smooth regression surface,
that is,
%
%
\begin{equation}
\label{eqregprob} y_i = \mu(x_i) + \varepsilon_i,
\qquad\varepsilon_i \sim\mbox{N}\bigl(0, \sigma^2\bigr).
\end{equation}
As motivated before, the true regression surface $\mu_0$ might depend
only on a subset of variables $I$ and have anisotropic smoothness in
the remaining variables. Accordingly, we assume $\mu_0 \in C^{\bfala
}[0,1]^I$ for some $I$ and $\alpha$. Also, assume the true value
$\sigma_0$ of $\sigma$ lies in an interval $[a, b] \subset[0, \infty)$.

We use the law of $W^{\bfaca}$ as a prior on $\mu$, with the prior on
$\mathbf{A}$ as in \textup{(P1)--(P4)}.
We also assume a prior on $\sigma$ supported on $[a, b]$. Denote the
posterior distribution by $\Pi( \cdot\mid y_1, \ldots, y_n)$. Let
$\llVert \mu\rrVert _n^2 = n^{-1} \sum_{i=1}^n \mu^2(x_i)$
denote the $L_2$ norm corresponding to the empirical distribution of
the design points. The posterior is said to contract at a rate
$\varepsilon_n$, if for every sufficiently large $M$,
%
%
\begin{equation}
\label{eqpostconcdefn} \mbox{E}_{\mu_0, \sigma_0} \Pi \bigl[(\mu, \sigma)\dvtx  \llVert \mu-
\mu_0\rrVert _n + \llvert \sigma- \sigma_0
\rrvert > M \varepsilon_n \mid y_1, \ldots, y_n
\bigr] \to0.
\end{equation}

%
%
\begin{teo}\label{teoregression}
Consider the nonparametric regression model (\ref{eqregprob}) with
$W^{\bfaca}$ as a prior on $\mu$, with the prior on $\mathbf{A}$ as
in \textup{(P1)--(P4)}. Let $\bfal= (\alpha_1, \ldots, \alpha
_d)$ be a vector of positive numbers and $I$ be a subset of $\{1,
\ldots, d\}$. If $\mu_0 \in C^{\bfala}[0,1]^I$, then~(\ref{eqpostconcdefn})
holds with $\varepsilon_n = n^{-\alpha_{0I}/(2
\alpha_{0I} + 1)} \log^{\kappa} n$, where $\alpha_{0I}^{-1} =\break
 \sum_{j \in I} \alpha_j^{-1}$ and \mbox{$\kappa> 0$} is a positive constant.
\end{teo}
Thus, one obtains the minimax optimal rate up to a log factor adapting
to the unknown dimensionality and anisotropic smoothness.


%
%
\begin{rem}
Theorem~\ref{teoregression} guarantees adaptive estimation of the
regression surface. It is often of interest additionally to select the
important variables affecting the response $y$. In the context of
variable selection in a normal linear model, \citet{barbieri2004optimal} advocated using the median probability model
consisting of the variables with posterior marginal inclusion
probability greater than or equal to half. They also proved the
predictive optimality of such models. The same approach could be
followed here for variable selection; however, the issue of optimality
needs to be studied.
\end{rem}

A similar result on adaptation holds for density estimation using the
logistic Gaussian process, where an unknown density $g$ on the
hypercube $[0, 1]^d$ is modeled as
%
%
\begin{equation}
\label{eqlogisticgauss} g(t) = \frac{e^{\mu(t)}}{{\int_{[0, 1]^d}} e^{\mu(s)}\,ds }
\end{equation}
for a function $\mu\dvtx  \mathbb{R}^d \to\mathbb{R}$.
Suppose $X_1, \ldots, X_n$ are drawn i.i.d. from a continuous,
everywhere positive density $g_0 = \log\mu_0$ on $[0, 1]^d$.\vadjust{\goodbreak}
%
%
\begin{teo}\label{teologistic}
Suppose one uses the law of $W^{\bfaca}$ as a prior on $\mu$ in (\ref
{eqlogisticgauss}), with the prior on $\mathbf{A}$ as in
\textup{(P1)--(P4)}. Let $\bfal= (\alpha_1, \ldots, \alpha_d)$ be a
vector of positive numbers and $I$ be a subset of $\{1, \ldots, d\}$.
If $\mu_0 = \log g_0 \in C^{\bfala}[0,1]^I$, then the posterior
contracts at the rate $\varepsilon_n = n^{-\alpha_{0I}/(2 \alpha_{0I} +
1)} \log^{\kappa} n$ with respect to the Hellinger distance, where
$\alpha_{0I}^{-1} = \sum_{j \in I} \alpha_j^{-1}$.
\end{teo}
The proofs of Theorems~\ref{teoregression} and~\ref{teologistic}
follow in a straightforward manner from our main results in Theorems
\ref{teoanisomain} and~\ref{teodimredmain}. We do not provide a
proof here since the steps are very similar to those in Section~3 of \citet{van2008rates}.

\subsection{Lower bounds on posterior contraction rates}\label{seclowerbounds}

In this section, we derive a lower bound to the posterior convergence
rate when the true function has anisotropic smoothness and a single
bandwidth GP is used, exhibiting the necessity of the multi-bandwidth
process. Consider once again the regression function estimation setting
(\ref{eqregprob}); we assume $\sigma$ to be fixed and known here.
For technical simplicity, we will formulate our lower bounds results in
a slightly different setting;\vspace*{1pt} we assume a random design with $x_i \sim
Q$, where $Q$ is a distribution on $[0, 1]^d$ admitting a density
$\omega$ with respect to the Lebesgue measure on $[0, 1]^d$. Let
$\llVert \mu\rrVert _{2, \omega}^2 = \int_{[0, 1]^d} \mu
^2(t) \omega(t) \,dt$. Recall that $\xi_n$ is a lower bound
[\citet{castillo2008lower}] to the posterior convergence rate
around $\mu_0$ in $\llVert \cdot\rrVert _{2, \omega}$ if
%
%
\begin{equation}
\label{eqlbdefn} P\bigl( \llVert \mu- \mu_0\rrVert _{2,\omega}
\leq\xi_n \mid Y_1, \ldots, Y_n,
x_1, \ldots, x_n\bigr) \to0 \qquad\mbox{a.s.}
\end{equation}
as $n \to\infty$.

In the following, for a positive random variable $A$ stochastically
independent of~$W$, we shall consider a rescaled Gaussian process
$W^{A}$ as the prior on the regression function $\mu$, with $A$
assigned the following prior (LBP):
\begin{longlist}[(LBP)]
\item[(LBP)]
Consider $A^{d} \sim h$, where $h$ is a density on the positive real
line satisfying for any $\varepsilon_1, \varepsilon_2 >0$ sufficiently small,
%
%
\begin{eqnarray}
h(x) &\leq& B_0 \exp(-B_1 /x), \qquad x <
\varepsilon_1,\label{eqgdens3.11}
\\
C_1 x^p \exp\bigl(-D_1 x
\log^q x\bigr) &\leq& h(x) \leq C_2 x^p \exp
\bigl(-D_2 x \log^q x\bigr),
\nonumber\\[-10pt] \label{eqgdens1}\\[-10pt]
\eqntext{x > 1/\varepsilon_2}
\end{eqnarray}
for positive constants $B_0, B_1, C_1, C_2, D_1, D_2$ and constants
$p,q \geq0$.
\end{longlist}
The tail behavior of $A$ implied by (\ref{eqgdens1}) is exactly the
same as equation (3.4) in \citet{van2009adaptive}. To establish the
lower bound result, we further need
to control the behavior of $h$ near zero, which is specified in (\ref
{eqgdens3.11}). In particular, letting $h$ to be a gamma density truncated
to $[c, \infty)$ for any constant $c > 0$ would satisfy (LBP).
Examples of $h$ supported on $(0, \infty)$ satisfying (LBP)
include\vadjust{\goodbreak}
the three-parameter generalized inverse Gaussian (gIG) family with
probability density function
\[
h(x) = \frac{(a/b)^{p/2}}{2 K_p(\sqrt{ab})} x^{p-1}e^{-(a x +
b/x)/2}, \qquad x > 0,
\]
where $K_p$ is a modified Bessel function of the second kind, $a > 0, b
> 0$ and $p$ is a
real parameter. Thus, being a subclass of prior distributions on $A$
considered by \citet{van2009adaptive}, (LBP) results in the minimax
rate of posterior contraction upto a logarithmic term adaptively over
isotropic $\alpha$-H\"{o}lder functions of $d$-variables for any
$\alpha> 0$ with respect to the empirical $L_2$ norm.

Next, consider a class of functions (LBT) for the true regression
function~$\mu_0$:
\begin{longlist}[(LBT)]
\item[(LBT)]
For integers $d \geq3$, $1 \leq d_1 < d$ and $d_2 = d - d_1$, define
$I_1 = \{1, \ldots, d_1\}$ and $I_2 = \{d_1 +1, \ldots, d\}$. Assume
$\mu_0(t) = \zeta(t_{I_1})\eta(t_{I_2})$,\vspace*{1pt}
where $\zeta\in C^{\alpha}[0,1]^{d_1}$, has support 
$[v,w]^{d_1}$ for $0 < v < w <1$ and $\eta\dvtx  [0, 1]^{d_2} \to\mathbb
{R}$ is infinitely differentiable with support $[v, w]^{d_2}$.
Assume further that there exists $0 < u < 1$, $\beta\geq\alpha$,
such that
for any constants $K_{1}$ and
$K_{2}$ sufficiently large,
the average (in $L_2$ sense) tail of the Fourier transform $\hat{\mu
}_0$ of $\mu_0$ satisfies
%
%
\begin{equation}
\label{eqmu0condnew} \int_{\llVert \lambda_{I_1}\rrVert > K_{1}, \llVert
\lambda_{I_2}\rrVert \in[K_2, 2K_2] } \bigl\llvert \hat{
\mu}_0(\lambda)\bigr\rrvert ^2 \,d\lambda\geq
K_{1}^{-2\beta} \exp\bigl\{-2 K_{2}^u
\bigr\}.
\end{equation}
\end{longlist}
Since $\hat{\mu}_0(\lambda) = \hat{\zeta}(\lambda_{I_1}) \hat
{\eta}(\lambda_{I_2})$, it suffices to have $\zeta\in C^{\alpha}[0,
1]^{d_1}$ and $\eta$\break infinitely  smooth\vspace*{1pt} such that $\int_{\llVert
\lambda_{1}\rrVert > K_{1}} | \hat{\zeta}(\lambda_{1}) |^2 \,d
\lambda_1 \geq K_{1}^{-2\beta}$ and\break  $\int_{\llVert \lambda
_{2}\rrVert \in[K_{2}, 2 K_2]} | \hat{\eta}(\lambda_{2}) |^2
\,d \lambda_2 \geq\exp\{-2 K_{2}^u \}$. We provide examples of
functions $\zeta$ and $\eta$ satisfying (LBT) in Appendix~\ref{app0}
for appropriate $\alpha, \beta$ and $u$; refer also to the discussion
in the last paragraph of this section.

In the following Theorem~\ref{teolb}, we demonstrate that (LBP) can
lead to a slower (in an
exponent of $n$) rate of convergence compared to the multi-bandwidth
case with respect to the $\llVert \cdot\rrVert _{2,\omega}$
norm when the true regression function $\mu_0$ belongs to the
anisotropic class (LBT).
%
%
\begin{teo}\label{teolb}
Consider the setting (\ref{eqregprob}) with (LBP) as the prior and
assume the true regression function $\mu_0$ satisfies (LBT). Let
$\gamma_1 = \alpha/(2 \alpha+ d)$, $\gamma_2 = \alpha/(2 \alpha+
d_1)$. Assume that the exponent $\beta$ appearing in (\ref
{eqmu0condnew}) satisfies $\beta< \alpha\gamma_2/ \gamma_1$. Then~(\ref{eqlbdefn}) holds with $\xi_n = n^{-\gamma}$ for any $\gamma
\in(\gamma_1, \gamma_2)$ satisfying $\beta< \alpha\gamma/\gamma_1$.
\end{teo}

Theorem~\ref{teoanisomain} implies\vspace*{1pt} that the posterior rate of
convergence for estimating $\mu_0$ satisfying (LBT) under the proposed
multi-bandwidth GP is faster than $n^{-\psi/(2\psi+ 1)}$ with $1/\psi
= d_1/\alpha+ d_2/m$ for any large integer $m$. Hence, this rate can
be made arbitrarily close to
$n^{-\gamma_2}$ (in an exponent of $n$). On the other hand, since $\mu
_0$ satisfying (LBT) belongs to the isotropic H\"{o}lder class
$C^{\alpha}[0, 1]^d$,\vadjust{\goodbreak} the posterior convergence rate under (LBP) is
bounded above by $n^{-\gamma_1}$ up to a logarithmic term [\citet
{van2009adaptive}].\footnote{We note once again that the upper bound
results are for the empirical $L_2$ norm.} Theorem~\ref{teolb} shows
that under suitable conditions, one can find $\gamma\in(\gamma_1,
\gamma_2)$
such that $n^{-\gamma}$ is a lower bound (\ref{eqlbdefn}) to the
posterior convergence rate under (LBP) in the $\llVert \cdot\rrVert _{2, \omega}$ norm. In other words, the posterior is
concentrated in the annulus with outer radius $n^{-\gamma_1}$ and
inner radius $n^{-\gamma}$. This exhibits the lack of efficiency
incurred by using a single bandwidth,\footnote{Albeit in a slightly
different norm.} since $n^{-\gamma}$ is slower than $n^{-\gamma_2}$ by
a genuine power of $n$.


The exponent $\beta$ in the $L_2$ averaged Fourier tail of $\zeta$
dictates the lower bound~$\xi_n$ in Theorem~\ref{teolb}, and hence
the assumption $\beta< \alpha\gamma_2/\gamma_1$ merits a
discussion. Fix~$\alpha, d_1$. Since $\gamma_2/\gamma_1 = 1 + d_2/(2
\alpha+ d_1)$, one can allow $\beta= K \alpha$ for any positive
constant $K > 1$ by choosing $d_2$ large enough in Theorem~\ref{teolb}. Intuitively, more the number of dimensions that contribute to
the anisotropy, less stringent becomes the requirement on $\zeta$.
In Appendix~\ref{app0}, we present a specific example of $\zeta$ in
the case $d_1 =1$ with $\alpha=1$ and $\beta= 3/2$. In this case,
$\beta< \alpha\gamma_2/ \gamma_1$ is equivalent to $3/2 < 1 +
d_2/3$, which is satisfied for any $d_2 \geq2$.
\section{Auxiliary results}
In this section, we present a number of auxiliary results that are
crucially used to prove the main results.

\subsection{Properties of the multi-bandwidth Gaussian process}\label{secproprescaled}
We first summarize some properties of the RKHS of the scaled process
$W^{\bfaa}$ for a fixed vector of scales~$\mathbf{a}$. Lemmas~\ref{lemrkhschar}--\ref{lemscaledentropy} generalize the results in
Section~4 of \citet{van2009adaptive} from a single scaling to a
vector of scales; we briefly sketch the proofs emphasizing the places
we differ substantially from \citet{van2009adaptive}.


Assume that the spectral measure $\nu$ of $W$ has a spectral density $f$.
For $\mathbf{a}\in\mathbb{R}_{*}^d$, the rescaled process $W^{\bfaa
}$ has a spectral measure $\nu_{\mathbf{a}}$ given by $\nu_{\mathbf
{a}}(B) = \nu(B./{\mathbf{a}})$. Further, $\nu_{\mathbf{a}}$ admits
a spectral density $f_{\mathbf{a}}$, with\vspace*{1pt} $f_{\mathbf{a}}(\lambda) =
\frac{1}{\mathbf{a}^*} f(\lambda./\mathbf{a})$. For $w_0 \in C[0,
1]^d$, define $\phi^{\mathbf{a}}_{w_0}(\varepsilon)$ to be the
concentration function of the rescaled Gaussian process~$W^{\bfaa}$.

As a straightforward extension of Lemmas 4.1 and 4.2 in \citet
{van2009adaptive}, it turns out that the RKHS of the process $W^{\bfaa
}$ can be characterized as below.
%
%
\begin{lem}\label{lemrkhschar}
The RKHS $\mathbb{H}^{\mathbf{a}}$ of the process $\{W_t^{\mathbf
{a}}\dvtx  t \in[0, 1]^d\}$ consists of real parts of the functions
\[
t \mapsto\int e^{ i (\lambda, t)} g(\lambda) \nu_{\mathbf
{a}}(d\lambda),
\]
where $g$ runs over the complex Hilbert space $L_2(\nu_{\mathbf
{a}})$. Further, the RKHS norm of the element in the above display is
given by $\llVert  g\rrVert _{L_2(\nu_{\mathbf{a}})}$.\vadjust{\goodbreak}
\end{lem}

Lemma 4.3 of \citet{van2009adaptive} shows that for any isotropic
H\"{o}lder smooth function $w$, convolutions with an appropriately
chosen class of higher order kernels indexed by the scaling parameter
$a$ belong to the RKHS. This suggests that driving the bandwidth $1/a$
to zero, one can obtain improved approximations to any H\"{o}lder
smooth function. The following Lemma~\ref{lemscaledapprox}
illustrates the usefulness of using separate bandwidths for each
dimension for approximating anisotropic H\"{o}lder functions from the RKHS.

%
%
\begin{lem}\label{lemscaledapprox}
Assume $\nu$ has a density with respect to the Lebesgue measure which
is bounded away from zero on a neighborhood of the origin. Let $\bfal
\in\mathbb{R}_*^d$ be given. Then, for any subset $I$ of $\{1, \ldots, d\}$ and $w \in C^{\bfala}[0, 1]^I$, there exists constants $C$ and
$D$ depending only on $\nu$ and $w$ such that, for $a_i$'s large enough,
\[
\inf \biggl\{ \llVert h\rrVert _{\mathbb{H}^{\mathbf{a}}}^2\dvtx  \llVert h-w
\rrVert _{\infty} \leq C \sum_{i \in I}
a_i^{-\alpha_i} \biggr\} \leq D \mathbf{a}^*.
\]
\end{lem}
\begin{pf}
We shall prove the result for $w \in C^{\bfala}[0, 1]^d$ and sketch an
argument for extending the proof to any $w \in C^{\bfala}[0, 1]^I$.

Let $\psi_j, j = 1, \ldots, d$, be a set of higher order kernels
which satisfy\break  $\int\psi_j(t_j) \,dt_j = 1$, $\int t_j^k \psi_j(t_j)
\,dt_j = 0$ for any positive integer $k$ and\break  $\int|t_j|^{\alpha_j}
|\psi_j(t_j)| \,dt_j \leq1$. Define\vspace*{1pt} $\psi\dvtx  \mathbb{R}^d \to\mathbb
{C}$ by $\psi(t) = \psi_1(t_1)\cdots\psi_d(t_d)$ so that one has
$\int_{\mathbb{R}^d} \psi(t) \,dt = 1$, $\int_{\mathbb{R}^d} t^{k}
\psi(t) \,dt = 0$ for any nonzero multiindex $k = (k_1, \ldots, k_d)$,
and\vspace*{1pt} the functions $|\hat{\psi}|/f$ and $|\hat{\psi}|^2/f$ are
uniformly bounded, where $\hat{\psi}$ denotes the Fourier transform
of $\psi$.

For a vector of positive numbers $\mathbf{a}= (a_1, \ldots, a_d)$,
let $\psi_{\mathbf{a}}(t) = \mathbf{a}^*   \psi(\mathbf{a}\cdot
t)$, where $\mathbf{a}^* = \prod_{j=1}^d a_j$.
Proceeding as in Lemma 4.3 of \citet{van2009adaptive}, one can
show that the convolution $\psi_{\mathbf{a}}*w$ is contained in the
RKHS $\mathbb{H}^{\mathbf{a}}$ and the squared RKHS norm of $\psi
_{\mathbf{a}}*w$ is bounded by $D \mathbf{a}^*$, with $D$ depending
only on $\nu$ and $w$.
Thus, the proof of Lemma~\ref{lemscaledapprox} would be completed if
we can show that
$\llVert \psi_{\mathbf{a}}*w - w\rrVert _{\infty} \leq C
\sum_{j=1}^d a_j^{-\alpha_j}$.

We have, for any $t \in\mathbb{R}^d$,
\[
\psi_{\mathbf{a}}*w(t) - w(t) = \int\psi(s) \bigl\{w(t - s./\mathbf{a}) - w(t)
\bigr\} \,ds.
\]
For $1 \leq j \leq d-1$, let $u^{(j)}$ denote\vspace*{1pt} the vector in $\mathbb
{R}^d$ with $u^{(j)}_i = 0$ for $i = 1, \ldots, j$ and $u^{(j)}_i = 1$
for $i = j+1, \ldots, d$. For any two vectors $x, y \in\mathbb
{R}^d$, we can navigate from $x$ to $y$ in a piecewise linear fashion
traveling parallel to one of the coordinate axes at a time. The
vertices of the path will be given by $x^{(0)} = x$, $x^{(j)} = u^{(j)}
\cdot x + (1 - u^{(j)}) \cdot y$ for $j = 1, \ldots, d-1$ and $x^{(d)}
= y$.

A multivariate Taylor expansion of $w(t - s./a)$ around $w(t)$ cannot
take advantage of the anisotropic smoothness of $w$ across different\vadjust{\goodbreak}
coordinate axes. Letting $x = t, y = t - s./a$ and $x^{(j)}, j= 0, 1,
\ldots, d$ as above, let us write $w(y) - w(x)$ in the following
telescoping form:
\begin{eqnarray*}
w(y) - w(x) &=& \sum_{j=1}^d w
\bigl(x^{(j)}\bigr) - w\bigl(x^{(j-1)}\bigr)
\\
&=& \sum
_{j=1}^d w_j\bigl( t_j -
s_j/a_j \mid x^{(j)} \bigr) - w_j
\bigl( t_j \mid x^{(j)} \bigr),
\end{eqnarray*}
where the functions $w_j$ are as defined in Section~\ref{secnotn},
with $w_j(t \mid x) = w(x_1,\ldots,\break  x_{j-1},t, x_{j+1}, \ldots, x_d)$ for
any $t \in\mathbb{R}$ and $x \in\mathbb{R}^d$.

Thus,
\[
w(t - s./a) - w(t) = \sum_{j=1}^d
\Biggl[ \sum_{i=1}^{\lfloor\alpha
_j \rfloor} D^i
w_j\bigl(t_j \mid x^{(j)} \bigr)
\frac{(-s_j/a_j)^i}{i !} + S_j(t_j, -s_j/a_j)
\Biggr],
\]
where $\llvert  S_j(t_j, -s_j/a_j)\rrvert \leq K s_j^{\alpha
_j} a_j^{-\alpha_j}$ by (\ref{eqanisocond2}),
for a constant $K$ depending on $\nu$ and $w$ but not on $t$ and $s$.
Combining the above, we have
\begin{eqnarray*}
\biggl\llvert \int\psi(s) \bigl\{w(t - s./\mathbf{a}) - w(t)\bigr\} \,ds\biggr
\rrvert &=& \Biggl\llvert \sum_{j=1}^d \int
\psi(s_j) S_j(t_j, -s_j/a_j)
\,ds_j\Biggr\rrvert
\\
&\leq& C \sum_{j=1}^d a_j^{-\alpha_j}.
\end{eqnarray*}

If, $w \in C^{\bfala}[0, 1]^I$ for some subset $I$ of $\{1, \ldots,
d\}$ with $\llvert  I\rrvert = \tilde{d}$, so that $w(t) =
w_0(t_I)$ for some $w_0 \in C^{\bfal_I}[0, 1]^{\tilde{d}}$, then the
conclusion follows trivially follows from the observation $\psi
_{\mathbf{a}} * w = \psi_{\mathbf{a}_I} * w_0$.
\end{pf}

We next study the centered small ball probability of the rescaled
process and the metric entropy of the unit RKHS ball.

%
%
\begin{lem}\label{lemscaledsmallball}
For any $a_0$ positive, there exists constants $C$ and $\varepsilon_0 >
0$ such that for $\mathbf{a}\geq a_0$ and $\varepsilon< \varepsilon_0$,
\[
-\log\mathrm{P} \bigl( \bigl\llVert W^{\bfaa}\bigr\rrVert _{\infty}
\leq\varepsilon \bigr) \leq C \mathbf{a}^* \biggl( \log\frac{\bar
{\mathbf{a}}}{\varepsilon}
\biggr)^{d+1}.
\]
\end{lem}

\begin{pf}
This follows from Theorem 2 in \citet{kuelbs1993metric} and Lemma 4.6
in \citet{van2009adaptive}. Proceeding as in Lemma 4.6 in \citet{van2009adaptive} and Lemma~\ref{lemscaledentropy}, we obtain
%
%
\begin{equation}
\label{eqsupbound} \phi^{\mathbf{a}}(\varepsilon) + \log0.5 \leq K_1
\mathbf{a}^* \biggl( \log\frac{\phi_{0}^{\mathbf{a}}(\varepsilon)}{\varepsilon} \biggr)^{1+d}
\end{equation}
for some constant $K_1 > 0$.
Note that with $L= [0, a_1]\times\cdots\times[0, a_d]$,
\begin{eqnarray*}
\phi_{0}^{\mathbf{a}}(\varepsilon) &=& -\log P\bigl(\bigl\llVert
W^{\bfaa
}\bigr\rrVert _{\infty} \leq\varepsilon\bigr) = -\log P\Bigl(
\sup_{t \in L}\llvert W_t\rrvert \leq\varepsilon\Bigr)
\\
&\leq&-\log P\Bigl(\sup_{t \in[0, \bar{\mathbf{a}}]^d}\llvert W_t\rrvert
\leq\varepsilon\Bigr) \leq K_2 \biggl(\frac{\bar{\mathbf
{a}}}{\varepsilon}
\biggr)^{\tau}
\end{eqnarray*}
for some constant $K_2$ and $\tau> 0$, where the last inequality
follows from the proof of Lemma 4.6 in \citet{van2009adaptive}.
Inserting this bound in~(\ref{eqsupbound}), we obtain the desired result.
\end{pf}


Let $\mathbb{H}_1^{\mathbf{a}} $ denote the unit ball in the RKHS of
$W^{\bfaa}$.
%
%
\begin{lem}\label{lemscaledentropy}
There exists a constant $K$, depending only on $\nu$ and $d$, such
that, for $\varepsilon< 1/2$,
\[
\log\mbox{N}\bigl(\varepsilon, \mathbb{H}_1^{\mathbf{a}}, \llVert
\cdot\rrVert _{\infty}\bigr) \leq K {\mathbf{a}}^* \biggl( \log
\frac{1}{\varepsilon} \biggr)^{d+1}.
\]
\end{lem}

\begin{pf}
By Lemma~\ref{lemrkhschar}, an element of $\mathbb{H}_1^{\mathbf
{a}}$ can be written as the real part of the function $h\dvtx  [0, 1]^d \to
\mathbb{C}$ given by
%
%
\begin{equation}
h(t) = \int e^{i(\lambda, t)}\psi(\lambda)\nu_{\mathbf
{a}}(d\lambda)
\end{equation}
for $\psi\dvtx  \mathbb{R}^{d} \to\mathbb{C}$ a function with $\int
\llvert \psi(\lambda)\rrvert ^2 \nu_{\mathbf{a}}(d\lambda
) \leq1$.

For $z \in\mathbb{C}^d$, continue to denote the function $z \mapsto
\int e^{(\lambda, z)}\psi(\lambda)\nu_{a}(d\lambda)$ by $h$. Using
the Cauchy--Schwarz inequality and the change of variable theorem,
%
%
\begin{equation}
\label{eqanaext} \bigl\llvert h(z)\bigr\rrvert ^2 \leq\int
e^{ (\lambda, 2 \mathbf
{a}\cdot\operatorname{Re}(z)) } \nu(d\lambda),
\end{equation}
where $\operatorname{Re}(z)$ denotes the vector whose $j$th element is the
real part of $z_j$ for $j = 1, \ldots, d$, and $\mathbf{a}\cdot\operatorname{Re}(z) = (a_1 \operatorname{Re}(z_1), \ldots, a_d \operatorname{Re}(z_d))^{\T}$.
From (\ref{eqanaext}) and the dominated convergence theorem, any $h
\in\mathbb{H}_1^{\mathbf{a}}$ can be analytically extended to
$\Gamma= \{ z \in\mathbb{C}^d\dvtx  \llVert 2 \mathbf{a}\cdot\operatorname{Re}(z)\rrVert _2 < \delta\}$. Clearly, $\Gamma$ contains a
strip $\Omega$ in $\mathbb{C}^d$ given by $\Omega= \{z \in\mathbb
{C}^d\dvtx  \llvert \operatorname{Re}(z_j)\rrvert \leq R_j, j =1, \ldots, d\}$ with $R_j = \delta/(6a_j \sqrt{d})$. Also, for every $z \in
\Omega$, $h$ satisfies the uniform bound $\llvert  h(z)\rrvert
^2 \leq\int e^{\delta\llVert \lambda\rrVert } \nu(d
\lambda) = C^2$.


Let $R = (R_1, \ldots, R_d)^{\T}$. Partition $T = [0, 1]^d$ into
rectangles $\Gamma_1, \ldots, \Gamma_m$ with centers $\{t_1, t_2,
\ldots, t_m\}$ such that given any $z \in T$, there exists $\Gamma_j$
with center $t_j = (t_{j1}, \ldots, t_{jd})^{\T}$ with $\llvert
z_i - t_{ji}\rrvert \leq R_i/4, i = 1, \ldots, d$.
Consider the piecewise polynomials $P = \sum_{j=1}^{m} P_{j, \gamma
_j} 1_{\Gamma_j}$ with
\[
P_{j, \gamma_j}(t) = \sum_{n. \leq k}
\gamma_{j,n} (t-t_j)^{n}.
\]
A finite set of functions $\mathcal{P}_{\mathbf{a}}$ is obtained by
discretizing the coefficients $\gamma_{j,n}$ for each $j$ and $n$ over
a grid of mesh width $\varepsilon/ R^{n}$
in the interval $[-C/R^{n}, C/R^{n}]$, with $R^n = R_1^{n_1} \cdots
R_d^{n_d}$ and $C$ defined as above. Choosing $m \precsim1/ R^*$ and
$k$ such that $k \precsim\log(1/\varepsilon)$ and $(2/3)^k \leq K
\varepsilon$ for some constant $K$, the collection $\mathcal{P}_{\mathbf
{a}}$ can be shown to form a $K \varepsilon$-net to $\mathbb
{H}_1^{\mathbf{a}}$ using (\ref{eqtailest}) and (\ref{eqpolyapprox}):
%
%
\begin{eqnarray}
\biggl\llvert \sum_{n.> k} \frac{D^nh_{\psi}(t_i)}{n!}(z -
t_i)^n \biggr\rrvert &\leq&\sum
_{n. > k} \frac{C}{R^n} (R/2)^n \leq KC \biggl(
\frac{2}{3} \biggr)^k, \label{eqtailest}
\\
\biggl\llvert \sum_{n. \leq k} \frac{D^nh_{\psi}(t_i)}{n!}(z -
t_i)^n - P_{i,\gamma_i}(z) \biggr\rrvert &\leq& K
\varepsilon. \label{eqpolyapprox}
\end{eqnarray}
The details are similar to Lemma 4.5 in \citet{van2009adaptive}, and
hence omitted.
\end{pf}

Lemma 4.7 of \citet{van2009adaptive} exploited a containment relation
among unit RKHS balls with different scalings to construct the sieves~$B_n$. Such a result sufficed in the single bandwidth case exploiting
the ordering of~$\mathbb{R}_+$. However, the result can only be
generalized with respect to the partial order on~$\mathbb{R}_+^d$~and
one does not obtain a straightforward generalization of their sieve
construction in the multi-bandwidth case since the entropy of their
sieve blows up in trying to control the joint probability of the
rescaling vector $\mathbf{a}$ outside a hyper-rectangle in~$\mathbb{R}_+^d$.

The problem mentioned above is fundamentally due to the curse of
dimensionality and one needs a more careful construction of the sieve
to avoid this problem.
In the proof of Lemma~\ref{lemscaledentropy}, a collection of
piece-wise polynomials is used to cover the unit RKHS ball $\mathbb
{H}_1^{\mathbf{a}}$. The main idea in the next Lemma~\ref{lemunionentropy2} is to exploit the fact that the same set of
piecewise polynomials can also be used to cover $\mathbb{H}_1^{\mathbf
{b}}$ for $\mathbf{b}$ sufficiently close to $\mathbf{a}$. We then
come up with a careful choice of a compact subset $\mathcal{Q}$ of
$\mathbb{R}_+^d$ that balances the metric entropy of the collection of
unit RKHS balls $\mathbb{H}_1^{\mathbf{a}}$ with $\mathbf{a}\in
\mathcal{Q}$ and the complement probability of $\mathcal{Q}$ under
the joint prior on $\mathbf{a}$.

%
%
\begin{figure}

\includegraphics{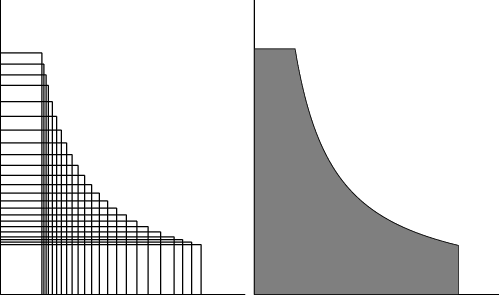}

\caption{Left panel: for $d=2$ and fixed $r > 1$, rectangles
$\mathcal{C}_{r^\theta} = \{ 0 \leq\mathbf{a}\leq r^{\theta}\}$
for different values of $\theta\in\mathcal{S}_{d-1}$. Right panel:
the region $\mathcal{Q}$ (shaded) resulting from the union of all such
rectangles.}\label{fig1}
\end{figure}

For $\bfu\in\mathbb{R}_+^d$, let $\mathcal{C}_{\bfu}$ denote the
rectangle in the positive quadrant given by $\mathbf{a}\leq\bfu$,
that is, $0 \leq a_j \leq u_j$ for all $j =1, \ldots, d$. For a fixed
$r > 1$, let $\mathcal{Q}^{(r)}$ consist of vectors $\mathbf{a}$ with
$a_j \leq r^{\theta_j}$ for some $\theta\in\mathcal{S}_{d-1}$. It
is easy to see that $\mathcal{Q}^{(r)}$ is a union of rectangles
$\mathcal{C}_{r^{\theta}}$ with $\theta$ varying over $\mathcal{S}_{d-1}$,
\[
\mathcal{Q}^{(r)} = \bigcup_{\theta\in\mathcal{S}_{d-1}} \mathcal
{C}_{r^{\theta}}.
\]
The outer boundary of $\mathcal{Q}^{(r)}$ consists of points $\mathbf
{a}$ with $0 < a_j \leq r$ for all $j = 1, \ldots, d$ and $\mathbf
{a}^* = r$ (see Figure~\ref{fig1}). By Lemma~\ref{lemscaledentropy}, for any such $\mathbf{a}$ in the outer boundary
of~$\mathcal{Q}^{(r)}$, the metric entropy of $\mathbb{H}_1^{\mathbf
{a}}$ is bounded by a constant multiple of $r \log^{d+1} (1/\varepsilon
)$. Our techniques were motivated by our observation that the entropy
remains of the same order even if one considers a union over the outer
boundary of $\mathcal{Q}^{(r)}$. We present a stronger result in Lemma
\ref{lemunionentropy2} which further states that the entropy remains
of the same order even if the union is considered over all of $\mathcal
{Q}^{(r)}$.


%
%
\begin{lem}\label{lemunionentropy2}
Let $\nu$ satisfy (\ref{eqsubexp}). Fix $r \geq1$. Then there
exists a constant $K$ depending on $\nu$ and $d$ only, so that, for
$\varepsilon< 1/2$,
\[
\log\mbox{N} \biggl(\varepsilon, \bigcup_{\mathbf{a}\in\mathcal
{Q}^{(r)}}
\mathbb{H}_1^{\mathbf{a}}, \llVert \cdot\rrVert _{\infty}
\biggr) \leq K r \biggl( \log\frac{1}{\varepsilon} \biggr)^{d+1}.
\]
\end{lem}


\begin{pf}
Let $\mathcal{Q} = \mathcal{Q}^{(r)}$. Fix $\mathbf{a}\in\mathcal
{Q}$. Let $I$ denote the subset of $\{1, \ldots, d\}$ such that $a_j
\leq1$ for all $j \in I$ and $a_j > 1$ for all $j \notin I$. Let
$\Omega^{\mathbf{a}} = \{z \in\mathbb{C}^d\dvtx  \llvert \operatorname{Re}(z_j)\rrvert \leq R_j, j =1, \ldots, d\}$ with $R_j = \delta
/(6a_j \sqrt{d})$ if $j \notin I$ and $R_j =\break  \delta/(6 \sqrt{d})$ if
$j \in I$. Note that for $z \in\Omega^{\mathbf{a}}$,
%
%
\begin{eqnarray}
\label{eqanaext1} \bigl\llVert 2 \mathbf{a}\cdot\operatorname{Re}(z)\bigr\rrVert
^2 &=& \sum_{j=1}^d 4
a_j^2 \bigl\llvert \operatorname{Re}(z_j)\bigr\rrvert
^2
\nonumber\\[-8pt]\\[-8pt]
& \leq&\sum_{j \in I} 4 R_j^2 +
\sum_{j \notin I} 4 a_j^2
R_j^2 \leq \frac{\delta^2}{9}.\nonumber
\end{eqnarray}
Hence, $\llVert 2 \mathbf{a}\cdot\operatorname{Re}(z)\rrVert \leq
\delta/3$ for $z \in\Omega^{\mathbf{a}}$. Following the argument
after the display in (\ref{eqanaext}), it thus follows that any
function $h \in\mathbb{H}_1^{\mathbf{a}}$ has an analytic extension
to $\Omega^{\mathbf{a}}$. Let $\mathbf{b}\in\mathcal{Q}$ satisfy
$\max_j \llvert  a_j - b_j\rrvert \leq1$. We shall exhibit
that any $h \in\mathbb{H}_1^{\mathbf{b}}$ can also be extended
analytically to the same strip $\Omega^{\mathbf{a}}$ by showing that
$\llVert 2 \mathbf{b}\cdot\operatorname{Re}(z)\rrVert < \delta$
on~$\Omega^{\mathbf{a}}$.
To that end, for $z \in\Omega^{\mathbf{a}}$, first observe that
%
%
\begin{eqnarray}
\label{eqanaext2} \bigl\llVert 2 \mathbf{b}\cdot\operatorname{Re}(z)\bigr\rrVert & \leq&\bigl
\llVert 2 \mathbf{a}\cdot\operatorname{Re}(z)\bigr\rrVert + \bigl\llVert 2 (\mathbf{b}-
\mathbf{a}) \cdot\operatorname{Re}(z)\bigr\rrVert.
\end{eqnarray}
The first term in the right-hand side above is bounded by $\delta/3$
following (\ref{eqanaext1}). To tackle the second term, we use
$\llvert  b_j - a_j\rrvert \leq1 \leq a_j$ for $j \notin I$
to conclude
%
%
\begin{eqnarray}
\label{eqanaext3}  \bigl\llVert 2 (\mathbf{b}- \mathbf{a}) \cdot\operatorname{Re}(z)\bigr
\rrVert ^2 &=& \sum_{j=1}^d 4
(b_j - a_j)^2 \bigl\llvert \operatorname{Re}(z_j)\bigr\rrvert ^2
\nonumber\\[-8pt]\\[-8pt]
& \leq&\sum_{j \in I} 4 R_j^2 +
\sum_{j \notin I} 4 a_j^2
R_j^2 \leq \frac{\delta^2}{9}.\nonumber
\end{eqnarray}
Combining (\ref{eqanaext2}) and (\ref{eqanaext3}), $\llVert 2
\mathbf{b}\cdot\operatorname{Re}(z)\rrVert \leq2\delta/3$ on $\Omega
^{\mathbf{a}}$, proving our claim. Since the same tail estimate as in
(\ref{eqtailest}) works for any $h \in\mathbb{H}_1^{\mathbf{b}}$,
it follows from (\ref{eqpolyapprox}) that the set of functions
$\mathcal{P}_{\mathbf{a}}$ form a $K\varepsilon$-net for $\mathbb
{H}_1^{\mathbf{b}}$.

Let $\mathcal{A}$ be a set of points in $\mathcal{Q}$ such that for
any $\mathbf{b}\in\mathcal{Q}$, there exists $\mathbf{a}\in
\mathcal{A}$ such that $\max_j{\llvert  a_j - b_j\rrvert }
\leq1$. One can clearly find an $\mathcal{A}$ with $\llvert
\mathcal{A}\rrvert \leq r^d$. The proof is completed by
observing that $\bigcup_{\mathbf{a}\in\mathcal{A}} \mathcal
{P}_{\mathbf{a}}$ form a $K \varepsilon$ net for $\bigcup_{\mathbf{a}\in
\mathcal{Q}} \mathbb{H}_1^{\mathbf{a}}$.
\end{pf}

\subsection{Results for lower bound}\label{ssecauxlb}
We now state and prove Lemma~\ref{lemlbnorm} which enables us to
derive a lower bound to the concentration function $\phi^{a}(\varepsilon
)$ for a fixed $a > 0$. This lower bound coupled with the model
identifiability of (\ref{eqregprob}) results in a lower bound to the
posterior concentration rate.

Denote by $\mathbb{H}^{a}$ the reproducing kernel Hilbert space of the
Gaussian process $W^{a}$. The key to obtaining a lower bound to the
concentration function $\phi_{\mu_0}^{a}(\varepsilon)$ when $\mu_0$
has anisotropic smoothness (LBT) is to find a lower bound to $\inf_{h
\in\mathbb{H}^a\dvtx  \llVert  h - \mu_0\rrVert _2 \leq\varepsilon
} \llVert  h\rrVert _{\mathbb{H}^a}^2/2$.

Let $h^a_{\psi}$ be the real part of the function $\int e^{i(\lambda,
t)} \psi(\lambda) f_a(\lambda)\,d\lambda$, where $f_a(\lambda)$
denotes the spectral density corresponding to the spectral measure $\nu
_a$ of $W^a$, so that $f_a(\lambda) = a^{-d} f(\lambda/a)$ with
$f(\lambda) = e^{-\llVert \lambda\rrVert ^2/4} /(2^d \pi
^{d/2})$. From \citet{van2009adaptive}, the RKHS of $\mathbb{H}^a$
consists of functions $h^a_{\psi}$ for $\psi\in L_2(\nu_a)$.
%
%
\begin{lem}\label{lemlbnorm}
If $\mu_0$ satisfies (LBT) for some $\beta> 0$ and $0< u <1$, then
for some constant $K > 0$,
\begin{eqnarray*}
&& \inf_{h \in\mathbb{H}^a\dvtx  \llVert  h - \mu_0\rrVert _2 \leq
\varepsilon} \frac{ \llVert  h\rrVert _{\mathbb{H}^a}^2}{2}
\\
&&\qquad \geq
\cases{
\displaystyle C \varepsilon^2 a^d\exp\bigl\{ K\varepsilon^{-2/\beta}/\bigl(4a^2\bigr) \bigr\}, &\quad if $a > \varepsilon^{-(2-u)/(2\beta)}$,
\vspace*{3pt}\cr
\displaystyle C\varepsilon^2 a^d\exp\bigl\{ K\varepsilon^{-u/\beta} \bigr\}, & \quad if $a \leq \varepsilon^{-(2-u)/(2\beta)}$.}
\end{eqnarray*}
\end{lem}


%
\begin{pf}
Let $r$ be a function such that $r$ is equal to $1$ on the support of
$\mu_0$, has itself support inside $[0, 1]^d$, $\llvert  r(\lambda
)\rrvert \leq1$ and the Fourier transform $\llvert \hat
{r}(\lambda)\rrvert \leq e^{-K \llVert \lambda\rrVert
^u}$ for large $\llVert \lambda\rrVert $. For any $0 < u <
1$, such an $r$ exists; we provide a construction for $u = 1/2$ in
Proposition~\ref{remexistr} in Appendix~\ref{app1}. Let $\psi\in
L_2(\nu_a)$ be such that $ \llVert  h_{\psi}^a - \mu_0\rrVert _2 < \varepsilon$. By construction, $h_\psi^a r$ has support
inside $[0,1]^d$ and $\mu_0 r = \mu_0$, so that $\llVert  h_\psi
^a r- \mu_0\rrVert _{2, \mathbb{R}} \leq
\llVert  h_{\psi}^a - \mu_0\rrVert _{2}$, where $\llVert
\cdot\rrVert _{2, \mathbb{R}}$ is the norm of $L_2(\mathbb
{R}^d)$ and $\llVert \cdot\rrVert _{2}$ the norm of
$L_2[0,1]^d$. The function $h_\psi^a r$ has Fourier transform\footnote
{According to our convention, $\widehat{f * g} = (2 \pi)^d \hat{f}
\hat{g}$. We drop the constant $(2 \pi)^d$ for notational
simplicity.} $(\psi_1 f_a) * \hat{r}$, where $\psi_1(\lambda) =
\psi(-\lambda)$. Hence, by Parseval's
identity $\llVert (\psi_1 f_a) * \hat{r} - \hat{\mu}_0\rrVert _{2, \mathbb{R}} < \varepsilon$.
Defining $\chi_{K_1, K_2} = I_{\{\lambda\in\mathbb{R}^d\dvtx  \llVert \lambda_{I_1}\rrVert > K_1, \llVert \lambda
_{I_2}\rrVert \in[K_2, 2 K_2]\}}$, we have
\begin{eqnarray*}
\llVert \psi_1 f_a * \hat{r} \chi_{2 K_1, K_2}
\rrVert _{2,
\mathbb{R}} &\geq& \llVert \hat{\mu}_0
\chi_{2 K_1, K_2}\rrVert _{2, \mathbb{R}} - \varepsilon.
\end{eqnarray*}
We also have from (LBT),
\begin{eqnarray*}
\llVert \hat{\mu}_0\chi_{2 K_1, K_2}\rrVert
^2_{2, \mathbb
{R}} &\geq& C \biggl(\frac{1}{K_1}
\biggr)^{2\beta} \exp\bigl\{-2K_2^{u}\bigr\}.
\end{eqnarray*}

Using the inequalities in the previous two displays and Proposition
\ref{lemmodlem16} in Appendix~\ref{app2} with $2 K_1$ instead of
$K_1$, $A_1 = K_1$, $f = \psi_1 f_a $ and $g= \hat{r}$,
%
%
\begin{eqnarray}
\label{eqL1}
&& \llVert \psi_1 f_a\chi_{K_1, \cdot}\rrVert _{2, \mathbb
{R}} \int_{\llVert  t_{I_1}\rrVert \leq K_1} \bigl\llvert \hat{r}(t)
\bigr\rrvert \,dt
\nonumber\\[-8pt]\\[-8pt]
&&\qquad \geq C \biggl(\frac{1}{K_1} \biggr)^{\beta}
e^{- K_2^{u}} - \varepsilon
- \llVert \psi_1 f_a\rrVert _{2, \mathbb{R}} \int
_{\llVert  t_{I_1}\rrVert > K_1} \bigl\llvert \hat {r}(t)\bigr\rrvert \,dt.\nonumber
\end{eqnarray}
Since $\llVert  h_{\psi}^a\rrVert _{\mathbb{H}^a} = \llVert \psi\sqrt{f_a}\rrVert _{2, \mathbb{R}}$ and $f_a$ is
symmetric about zero,
%
%
\begin{equation}
\label{eqL3} \llVert \psi_1 f_a\rrVert
_{2, \mathbb{R}} \leq C a^{-d/2}\bigl\llVert h_{\psi}^a
\bigr\rrVert _{\mathbb{H}^a}.
\end{equation}
Also,
%
%
\begin{eqnarray}
\label{eqL5} \llVert \psi_1 f_a \chi_{K_1, \cdot}
\rrVert _{2, \mathbb
{R}}^2 &=& \int_{\llVert \lambda_{I_1}\rrVert > K_1} \bigl(
\psi _1^2 f_a\bigr) f_a \leq
\Bigl(\sup_{\llVert \lambda_{I_1}\rrVert > K_1}f_a \Bigr)\bigl\llVert
h_{\psi}^a\bigr\rrVert ^2_{\mathbb{H}^a}
\nonumber\\[-8pt]\\[-8pt]
&\leq& Ca^{-d} e^{-K_1^2/(4a^2)} \bigl\llVert h_{\psi}^a
\bigr\rrVert ^2_{\mathbb{H}^a}.\nonumber
\end{eqnarray}
Equations (\ref{eqL3}) and (\ref{eqL5}) imply
%
%
\begin{equation}
\label{eqL6} \bigl\llVert h_{\psi}^a\bigr\rrVert
_{\mathbb{H}^a} \geq\frac
{a^{d/2} (C  (1/K_1)^{\beta} e^{-K_2^{u}} - \varepsilon )}{ e^{-K_1^2/(8a^2)} \int_{\llVert  t_{I_1}\rrVert \leq K_1} \llvert \hat
{r}(t)\rrvert  \,dt +\int_{\llVert  t_{I_1}\rrVert > K_1}
\llvert \hat{r}(t)\rrvert  \,dt}.
\end{equation}
For large $K_1$, $\int_{ \llVert  t_{I_1}\rrVert >K_1} \llvert \hat{r}(t)\rrvert  \,dt
\leq\int_{ \llVert  t\rrVert >K_1} \llvert \hat{r}(t)\rrvert  \,dt \leq\int_{
\llVert  t\rrVert >K_1} e^{-K \llVert  t\rrVert ^u} \,dt =\break  C
\int_{K_1}^{\infty} r^{d-1} e^{-K r^u} \,dr \leq e^{-K K_1^u}$ and
$\int_{ \llVert  t_{I_1}\rrVert \leq K_1} \llvert \hat
{r}(t)\rrvert  \,dt \leq\int\llvert \hat{r}(t)\rrvert
\,dt \leq C$. We can thus bound the denominator in (\ref{eqL6}) from
above by
$C e^{-K_1^2/(8a^2)} + e^{-K K_1^u}$. For fixed $K_1$, $C
e^{-K_1^2/(8a^2)} < e^{-K K_1^u} \Leftrightarrow a \leq C'
K_1^{(2-u)/2}$. Hence, the denominator can be bounded above by $2 C
e^{-K_1^2/(8a^2)}$ or $2 e^{-K K_1^u}$ depending on whether $a$ is,
respectively, larger or smaller than $C' K_1^{(2-u)/2}$.

Fix\vspace*{1pt} $K_2$ to be a large number. Choose $K_1$ large enough depending on
$\varepsilon$ such that $C (1/K_1)^{\beta} e^{-cK_2^{u}} = 2\varepsilon$;
this implies $K_1 =
(C/\varepsilon )^{1/ \beta}$ for some \mbox{$C > 0$}. With this
choice, $K_1^{(2-u)/2} = C (1/\varepsilon)^{(2-u)/(2 \beta)}$, $e^{-K
K_1^u} = \exp\{ -K (1/\varepsilon)^{u/ \beta} \}$ and
$e^{-K_1^2/(8a^2)} = \exp\{ - K (1/\varepsilon)^{2/\beta}/(8a^2) \}$.
Substituting the sequence of bounds in~(\ref{eqL6}), we have the
desired result.
\end{pf}

\section{Proof of main results}\label{secmainproofs}
We shall only provide a detailed proof of Theorem~\ref{teoanisomain}
and sketch the main steps in the proof of Theorem~\ref{teodimredmain}.
\subsection{\texorpdfstring{Proof of Theorem \protect\ref{teoanisomain}}{Proof of Theorem 3.1}}
Let us begin by observing that
\begin{eqnarray*}
\mathrm{P} \bigl(\bigl\llVert W^{\bfaca} - w_0\bigr\rrVert
_{\infty} \leq2\varepsilon \bigr) &=& \int\mathrm{P}\bigl(\bigl\llVert
W^{\bfaa} - w_0\bigr\rrVert _{\infty} \leq2\varepsilon
\bigr) \pi_{\bfaca}(d \mathbf {a})
\\
&=& \int \biggl\{ \int\mathrm{P}\bigl(\bigl\llVert W^{\bfaa} -
w_0\bigr\rrVert _{\infty} \leq2\varepsilon\bigr) \pi(\mathbf{a}
\mid\theta) \,d\mathbf{a} \biggr\} \pi(\theta) \,d\theta.
\end{eqnarray*}
As in \citet{van2009adaptive}, we first derive bounds on the
noncentered small ball probability for a fixed rescaling $\mathbf
{a}$, and then integrate over the distribution of $\mathbf{a}$ to
derive the same for $W^{\bfaca}$.

Given $\mathbf{a}\in\mathbb{R}_*^d$, recall the definition of the
centered and noncentered concentration functions of the process
$W^{\bfaa}$,
%
%
\begin{eqnarray}
\label{eqnoncentconc} \phi_0^{\mathbf{a}}(\varepsilon) &=& - \log\mathrm{P}
\bigl(\bigl\llVert W^{\bfaa}\bigr\rrVert _{\infty} \leq\varepsilon
\bigr),
\nonumber\\[-8pt]\\[-8pt]
\phi_{w_0}^{\mathbf{a}}(\varepsilon) &=& \inf_{h \in\mathbb
{H}^{\mathbf{a}}\dvtx  \llVert  h-w_0\rrVert _{\infty} \leq
\varepsilon}
\llVert h\rrVert _{\mathbb{H}^{\mathbf{a}}}^2 - \log\mathrm{P}\bigl(\bigl\llVert
W^{\bfaa}\bigr\rrVert _{\infty} \leq \varepsilon\bigr).\nonumber
\end{eqnarray}
For a fixed $\mathbf{a}$, the noncentered small ball probability of
$W^{\bfaa}$ can be bound in terms of the concentration function as
follows [\citet{van2008reproducing}]:
\[
\mathrm{P}\bigl(\bigl\llVert W^{\bfaa} - w_0\bigr\rrVert
_{\infty} \leq 2\varepsilon\bigr) \geq e^{- \phi_{w_0}^{\mathbf{a}}(\varepsilon)}.
\]

Now, suppose that $w_0 \in C^{\bfala}[0, 1]^d$ for some $\bfal\in
\mathbb{R}_*^d$. From Lemmas~\ref{lemscaledapprox} and~\ref{lemscaledsmallball}, it follows that for every $a_0 > 0$, there
exist positive constants $\varepsilon_0 < 1/2$, $C, D$ and $E$ that
depend only on $w_0$ and $\nu$ such that, for $\mathbf{a}> a_0$,
$\varepsilon< \varepsilon_0$ and $C \sum_{i=1}^d a_i^{-\alpha_i} <
\varepsilon$,
\[
\phi_{w_0}^{\mathbf{a}}(\varepsilon) \leq D \mathbf{a}^* + E \mathbf
{a}^* \biggl( \log\frac{\bar{\mathbf{a}}}{\varepsilon} \biggr)^{1+d} \leq
K_1 \mathbf{a}^* \biggl( \log\frac{\bar{\mathbf{a}}}{\varepsilon
}
\biggr)^{1+d}
\]
with $K_1$ depending only on $a_0, \nu$ and $d$. Thus, for $\varepsilon<
\min\{\varepsilon_0, C_1 a_0^{-\bar{\alpha}}\}$, by (\ref
{eqnoncentconc}), for constants $K_2, \ldots, K_6 > 0$ and $C_2,
\dots, C_6 > 0$,
\begin{eqnarray*}
&&\mathrm{P} \bigl(\bigl\llVert W^{\bfaca} - w_0\bigr\rrVert
_{\infty} \leq2\varepsilon \bigr)
\\
&&\qquad \geq \int_{\theta} \biggl\{ \int e^{-\phi_{w_0}^{\mathbf{a}}(\varepsilon)} \pi(
\mathbf{a}\mid\theta) \,d\mathbf{a} \biggr\} \pi(\theta) \,d\theta
\\
&&\qquad \geq\int_{\theta} \biggl\{ \int_{a_1 = (C_1/\varepsilon)^{1/\alpha
_1}}^{2 (C_1/\varepsilon)^{1/\alpha_1}}
\cdots\int_{a_d =
(C_1/\varepsilon)^{1/\alpha_d}}^{2 (C_1/\varepsilon)^{1/\alpha_d}} e^{-K_1 \mathbf{a}^* \log^{1+d}(\bar{\mathbf{a}}/\varepsilon)} \pi (
\mathbf{a}\mid\theta) \,d\mathbf{a} \biggr\} \pi(\theta) \,d\theta
\\
&&\qquad \geq C_2 e^{- K_2 (1/\varepsilon)^{1/\alpha_0} \log^{1+d}(1/\varepsilon
) }
\\
&&\quad\qquad{}\times \int_{\theta} \biggl\{
\int_{a_1 = (C_1/\varepsilon)^{1/\alpha
_1}}^{2 (C_1/\varepsilon)^{1/\alpha_1}} \cdots\int_{a_d =
(C_1/\varepsilon)^{1/\alpha_d}}^{2 (C_1/\varepsilon)^{1/\alpha_d}}
\pi (\mathbf{a}\mid\theta) \,d\mathbf{a} \biggr\} \pi(\theta) \,d\theta.
\end{eqnarray*}
Let $\Gamma$ denote the region in the simplex $\mathcal{S}_{d-1}$
given by $\Gamma= \{\theta\in\mathcal{S}_{d-1}\dvtx    \tau< \theta_j
- \frac{\alpha_0}{\alpha_j} < 2 \tau, j = 1, \ldots, d-1\}$. Since
$\sum_{j=1}^d \alpha_0/\alpha_j = 1$, we can choose $\tau> 0$ small
enough to guarantee that any $\theta$ satisfying the set of
inequalities lies inside the simplex. Moreover, with $\theta_d = 1 -
\sum_{j=1}^{d-1} \theta_j$, one has $(d-1) \tau< \theta_d < 2(d-1)
\tau$. Choosing $\tau= C_3/\log(1/\varepsilon)$, one can show that
$\sum_{j=1}^d (1 / \varepsilon)^{1/(\alpha_j \theta_j)} \leq C_4
(1/\varepsilon)^{1/\alpha_0}$ for any $\theta\in\Gamma$.
Now,
\begin{eqnarray*}
&&\int \biggl\{ \int_{a_1 = (C_1/\varepsilon)^{1/\alpha_1}}^{2
(C_1/\varepsilon)^{1/\alpha_1}} \cdots\int
_{a_d = (C_1/\varepsilon
)^{1/\alpha_d}}^{2 (C_1/\varepsilon)^{1/\alpha_d}} \pi(\mathbf{a}\mid \theta) \,d\mathbf{a}
\biggr\} \pi(\theta) \,d\theta
\\
&&\qquad\geq\int \biggl\{ \int_{a_1 = (C_1/\varepsilon)^{1/\alpha_1}}^{2
(C_1/\varepsilon)^{1/\alpha_1}} \cdots\int
_{a_d = (C_1/\varepsilon
)^{1/\alpha_d}}^{2 (C_1/\varepsilon)^{1/\alpha_d}} e^{- \sum_{j=1}^d
a_j^{1/\theta_j}} \,d\mathbf{a} \biggr\}
\pi(\theta) \,d\theta
\\
&&\qquad\geq\int e^{- K_3 \sum_{j=1}^d (1/\varepsilon)^{1/\alpha_j \theta
_j}} \pi(\theta) \,d\theta
\\
&&\qquad \geq\int_{\theta\in\Gamma} e^{- K_4 (1/\varepsilon)^{1/\alpha_0}} \pi(\theta) \,d\theta\geq
C_5 e^{- K_5 (1/\varepsilon)^{1/\alpha_0}}.
\end{eqnarray*}
The last inequality in the above display uses that $\Gamma$ contains a
hyper-cube of width~$\tau$ so that its $\pi$-mass is at least
polynomial in $\tau$.
Hence,
%
%
\begin{equation}
\label{concprobfinal} \mathrm{P} \bigl(\bigl\llVert W^{\bfaca} - w_0
\bigr\rrVert _{\infty} \leq2\varepsilon \bigr) \geq C_6
e^{- K_6 (1/\varepsilon)^{1/\alpha_0}
\log^{1+d}(1/\varepsilon) }.
\end{equation}

Let $\mathbb{B}_1$ denote the unit sup-norm ball of $C[0, 1]^d$. For a
vector $\theta\in\mathcal{S}_{d-1}$ and positive constants $M, r >
1, \varepsilon$, let $B^{\theta} = B^{\theta}(M, r, \varepsilon)$ denote
the set,
\[
B^{\theta} = \bigcup_{\mathbf{a}\leq r^\theta} \bigl(M \mathbb
{H}_1^{\mathbf{a}}\bigr) + \varepsilon\mathbb{B}_1,
\]
where $r^{\theta}$ denotes the vector whose $j$th element is
$r^{\theta_j}$. We further let
\[
B = \bigcup_{\mathbf{a}\in\mathcal{Q}^{(r)}} \bigl(M \mathbb
{H}_1^{\mathbf{a}}\bigr) + \varepsilon\mathbb{B}_1 =
\bigcup_{\theta\in
\mathcal{S}_{d-1}} \bigcup_{\mathbf{a}\leq r^\theta}
\bigl(M \mathbb {H}_1^{\mathbf{a}}\bigr) + \varepsilon
\mathbb{B}_1.
\]
Let us first calculate the probability $\mathrm{P}(W^{\bfaca} \notin
B^{\theta} \mid\theta)$. Note that
\begin{eqnarray*}
\mathrm{P}\bigl(W^{\bfaca} \notin B^{\theta} \mid\theta\bigr) &=& \int
\mathrm{P}\bigl(W^{\theta} \notin B^{\theta} \bigr) \pi(\mathbf{a}\mid
\theta) \,d\mathbf{a}
\\
& \leq&\int_{\mathbf{a}\leq r^{\theta}} \mathrm{P}\bigl(W^{\bfaa} \notin
B^{\theta}\bigr) \pi(\mathbf{a}\mid\theta) \,d\mathbf{a}+ \mathrm{P}\bigl(
\mathbf{A}\nleq r^{\theta} \mid\theta\bigr),
\end{eqnarray*}
where $\mathrm{P}(W^{\bfaa} \nleq r \mid\theta)$ is a shorthand
notation for $\mathrm{P}(\mbox{at least one }   A_j > r^{\theta_j}
\mid\theta)$.

To tackle the first term in the last display, note that $B^{\theta}$
contains the set $M \mathbb{H}_1^{\mathbf{a}} + \varepsilon\mathbb
{B}_1$ for any $\mathbf{a}\leq r^{\theta}$. Hence, for any $\mathbf
{a}\leq r^{\theta}$, by Borell's inequality,
\begin{eqnarray*}
\mathrm{P}\bigl(W^{\bfaa} \notin B^{\theta}\bigr) &\leq&\mathrm{P}
\bigl(W^{\bfaa} \notin M\mathbb{H}_1^{\mathbf{a}} +
\varepsilon\mathbb{B}_1\bigr)
\\
& \leq& 1 - \Phi \bigl\{M + \Phi^{-1} \bigl(e^{-\phi_0^{\mathbf
{a}}(\varepsilon)} \bigr)
\bigr\}
\\
& \leq& 1 - \Phi \bigl\{M + \Phi^{-1} \bigl(e^{-\phi_0^{r^{\theta
}}(\varepsilon)} \bigr)
\bigr\} \leq e^{-\phi_0^{r^{\theta}}(\varepsilon)},
\end{eqnarray*}
if $M \geq- 2 \Phi^{-1} (e^{-\phi_0^{r^{\theta}}(\varepsilon)}
)$.
The penultimate inequality in the above display follows from the fact
that, with $T = [0, 1]^d$,
\[
e^{-\phi_0^{\mathbf{a}}(\varepsilon)} = \mathrm{P}\Bigl(\sup_{t \in\mathbf
{a}\cdot T} \llvert
W_t\rrvert \leq\varepsilon\Bigr) \geq\mathrm{P}\Bigl(\sup
_{t \in r^{\theta} \cdot T} \llvert W_t\rrvert \leq \varepsilon\Bigr) =
e^{-\phi_0^{r^{\theta}}(\varepsilon)}.
\]
By Lemma 4.10 of \citet{van2009adaptive}, $\Phi^{-1}(u) \geq\break  - \{2
\log(1/u)\}^{1/2}$ for $u \in(0,1)$. Hence, the last inequality in
the above display remains valid if we choose
\[
M \geq4 \sqrt{\phi_0^{r^{\theta}}(\varepsilon)}.
\]

Since $A_j^{1/\theta_j}$ follows a gamma distribution given $\theta
_j$, in view of Lemma 4.9 of \citet{van2009adaptive}, for $r$ larger
than a positive constant depending only on the parameters of the gamma
distribution,
\[
\mathrm{P}\bigl(A_j > r^{\theta_j} \mid\theta\bigr) \leq
C_1 r^{D_1} e^{- D_2 r}.
\]

Combining the above, since $B$ contains $B^{\theta}$ for every $\theta
\in\mathcal{S}_{d-1}$,
%
%
\begin{eqnarray}
\label{eqcomlprobfinal} \mathrm{P}\bigl(W^{\bfaca} \notin B\bigr) &=& \int
_{\theta} \biggl\{ \int\mathrm{P}\bigl(W^{\bfaa} \notin B
\mid\theta\bigr) g(\mathbf{a}\mid\theta) \biggr\}
\nonumber
\\
& \leq&\int_{\theta} \biggl\{ \int\mathrm{P}\bigl(W^{\bfaa}
\notin B^{\theta} \mid\theta\bigr) g(\mathbf{a}\mid\theta) \biggr\}
\\
& \leq& C_2 r^{D_1} e^{- D_2 r} + e^{-D_3 r \log(r/\varepsilon)^{d+1}}.\nonumber
\end{eqnarray}

From Lemma~\ref{lemunionentropy2}, the entropy of $B$ can be
estimated as
%
%
\begin{eqnarray}
\label{eqentropyfinal} \log N\bigl(2 \varepsilon, B, \llVert \cdot\rrVert _{\infty}
\bigr) & \leq& \log N\biggl(\varepsilon, \bigcup_{\theta\in\mathcal{S}_{d-1}}
\bigcup_{\mathbf{a}\leq r^\theta} \bigl(M \mathbb{H}_1^{\mathbf{a}}
\bigr), \llVert \cdot\rrVert _{\infty}\biggr)
\nonumber\\[-8pt]\\[-8pt]
& \leq& K r \log \biggl(\frac{M}{\varepsilon} \biggr)^{d+1}.\nonumber
\end{eqnarray}
Thus, (\ref{concprobfinal}), (\ref{eqcomlprobfinal}) and (\ref
{eqentropyfinal}) can be simultaneously satisfied if we choose, for
constants $\kappa, \kappa_1, \kappa_2 > 0$,
\begin{eqnarray*}
\varepsilon_n &=& n^{-\alpha_0/(2 \alpha_0 + 1)} \log^{\kappa}(n),
\\
r_n &=& n^{1/(2 \alpha_0 + 1)} \log^{\kappa_1}(n),
\\
M_n &=& r_n \log^{\kappa_2}(n).
\end{eqnarray*}

\subsection{\texorpdfstring{Proof of Theorem \protect\ref{teodimredmain}}{Proof of Theorem 3.2}}
For ease of notation, we shall make the simplifying assumption that the
random variable $B$ is degenerate at $1$. For $a > 0$ and $S \subset\{
1, \ldots, d\}$, let $\mathbb{H}^{a,S}$ denote the RKHS of $W^{\bfaa
}$, where $a_j = a$ for $j \in S$ and $a_j = 1$ for $j \notin S$.

For a subset $S \subset\{1, \ldots, d\}$ with $\llvert  S\rrvert = \tilde{d}$, and given positive constants $M, r, \xi,
\varepsilon$, let
\begin{eqnarray*}
B_{S} &=& B_S(M, r, \xi, \varepsilon)
\\
& =& \biggl[M \biggl(\frac{r}{\xi} \biggr)^{\tilde{d}}
\mathbb{H}_1^{r,
S} + \varepsilon\mathbb{B}_1
\biggr] \cup \biggl[\bigcup_{a < \xi} \bigl(M
\mathbb{H}_1^{a, S} \bigr) + \varepsilon\mathbb{B}_1
\biggr].
\end{eqnarray*}

Since, given $S$, $A^{\tilde{d}} \sim\operatorname{gamma}(b_1, b_2)$, it can
be shown that, for some constant $C_1 > 0$,
\[
\mathrm{P}\bigl(W^{\bfaca} \notin B_S \mid S\bigr) \precsim
e^{- C_1 r^{\tilde{d}^*}}.
\]
The dominating term in the $\varepsilon$ entropy of $B_S$ is bounded by
\[
C_2 r^{d^*} \log^{1+\tilde{d}} \biggl(\frac{C_3 M}{\varepsilon}
\biggr).
\]
While calculating the concentration probability around $w_0 \in
C^{\bfala}[0, 1]^I$, simply use the fact that $\operatorname{pr}(S = I) > 0$.

Combining the above, the sieves $B_n$ are constructed as
\[
B_n = \bigcup_{\tilde{d} = 1}^{d}
\bigcup_{S\dvtx  \llvert  S\rrvert = \tilde{d}} B_S\bigl(M_n^S,
r_n^S, \xi_n, \varepsilon_n\bigr),
\]
where, for constants $\kappa, \kappa_1 > 0$, $\varepsilon_n =
n^{-\alpha/(2 \alpha+ d_0)} \log^{\kappa} n, r_n^S =  (
n^{d_0/(2 \alpha+ d_0)}  )^{1/\llvert  S\rrvert }\*
\log^{\kappa_1}(n)$ and $(M_n^S)^2 = (r_n^S)^{\tilde{d}} \log
(r_n^S/\varepsilon_n)$.

\subsection{\texorpdfstring{Proof of Theorem \protect\ref{teolb}}{Proof of Theorem 3.7}}
Let $\gamma\in(\gamma_1, \gamma_2)$ satisfy $\beta< \alpha\gamma
/\gamma_1$.
We show that
%
%
\begin{equation}
\label{eqlbmain} \frac{P(\llVert  W^A - \mu_0\rrVert _2 < \xi_n)}{P(\llVert  W^A - \mu_0\rrVert _\infty< \varepsilon_n)} \precsim\exp\bigl\{ -2n\varepsilon_n^2
\bigr\}
\end{equation}
for $\xi_n = n^{-\gamma}$ and $\varepsilon_n = n^{-\alpha/(2\alpha
+d)}\log^{t_1}n$ for some appropriate constant $t_1 > 0$.
It then follows from the proof of Theorem 8 in \citet{van2011information} that $\xi_n$ is a lower bound to the rate of
posterior contraction around $\mu_0$.

We first derive an upper bound to $P(\llVert  W^A - \mu_0\rrVert _2 < \xi)$ for $\xi$ small. Let $g$ denote the density of $A$
induced from (LBP). Clearly, $P(\llVert  W^A - \mu_0\rrVert _2
< \xi) = \int_{a=0}^{\infty} P(\llVert  W^a - \mu_0\rrVert
_2 < \xi) g(a) \,da$.
First, find $u(\xi)$ sufficiently small and $v(\xi)$ sufficiently
large such that
%
%
\begin{eqnarray}
\label{lbeq1} P\bigl( A < u(\xi)\bigr) &\leq&\exp\bigl( -C \xi^{-d/\beta}
\bigr),\qquad
\nonumber\\[-8pt]\\[-8pt]
P\bigl( A > v(\xi)\bigr) &\leq& \exp\bigl\{ -C \xi^{-d/\beta}
\log^{s}(1/\xi)\bigr\}\nonumber
\end{eqnarray}
for some $s \in\mathbb{R}$.
From (LBP), a simple calculation yields
%
%
\begin{eqnarray}
\label{lbeq2}
P\bigl( A < u(\xi)\bigr) &\leq&\exp\bigl( -C u(\xi)^{-d}
\bigr),
\nonumber\\[-8pt]\\[-8pt]
P\bigl( A > v(\xi)\bigr) &\leq&\exp\bigl( -C v(\xi)^{d}
\log^q v(\xi)\bigr).\nonumber
\end{eqnarray}
Hence,\vspace*{-1pt} if we choose $u(\xi) = \xi^{1/\beta}$, and $v(\xi) = \frac
{\xi^{-1/\beta}}{2\sqrt{2}\log^{1/2}(1/\xi)}$, then (\ref
{lbeq1}) is satisfied with $s = q - d/2$.

For $a \in(u(\xi), v(\xi))$, we bound the noncentered small ball
probability $P(\llVert  W^a - \mu_0\rrVert _2 < \xi)$ above
by $\exp\{-\phi_{\mu_0}^a(\xi)\}$ [Lemma 5.3 of \citet{van2008reproducing}] and further invoke the lower bound to the
concentration function $\phi_{\mu_0}^a(\xi)$ developed in Lemma~\ref{lemlbnorm}.
Specifically, we subdivide $(u(\xi), v(\xi))$ into two disjoint
regions based on the conclusion of Lemma~\ref{lemlbnorm} with
$u=1/2$. If $a \in(u(\xi), \xi^{-3/(4\beta)})$, Lemma~\ref{lemlbnorm} implies
%
%
\begin{equation}
\label{lbeq3} \phi_{\mu_0}^a(\xi) \geq\xi^{2 + d/\beta}
\exp \bigl\{K \xi^{-1/(2\beta)} \bigr\} \geq\xi^{-d/\beta}.
\end{equation}
If $a \in( \xi^{-3/(4\beta)}, v(\xi))$,
then again from Lemma~\ref{lemlbnorm},
%
%
\begin{equation}
\label{lbeq4} \phi_{\mu_0}^a(\xi) \geq C
\xi^2a^d\exp\bigl\{ \xi^{-2/\beta}/
\bigl(4a^2\bigr)\bigr\}.
\end{equation}
%
%
Putting together all the bounds, and noting that $s \geq-d/2$,
%
%
\begin{eqnarray}\label{lbe56}
&& \int_{a=0}^{\infty}P\bigl( \bigl\llVert
W^a - \mu_0\bigr\rrVert _{2} < \xi
\bigr)g(a) \,da\nonumber
\\
&&\qquad \leq 2 \exp\bigl( - C \xi^{-d/\beta}\bigr)+ \exp\bigl\{- C
\xi^{-d/\beta} \log ^{-d/2}(1/\xi) \bigr\}
\\
&&\quad\qquad{}+ \int_{a=\xi^{-3/(4\beta)}}^{ v(\xi)} \exp\bigl[ -C\xi^2
a^d \exp\bigl\{ \xi^{-2/\beta}/\bigl(4a^2\bigr)\bigr
\} \bigr] g(a) \,da.\nonumber
\end{eqnarray}
Observe that $\psi(x) = x^d \exp\{\xi^{-2/\beta}/(4x^2)\}$ is
decreasing if $x \in(0, \xi^{-1/\beta}/ \sqrt{2d}]$ and
increasing if $x >\xi^{-1/\beta}/\sqrt{2d}$. Since $v(\xi) \leq\xi
^{-1/\beta}/ \sqrt{2d}$, the third term in the r.h.s. of (\ref
{lbe56}) is bounded above by $\exp [ - C\xi^2 \psi\{v(\xi)\}
]$. Since $\psi\{ v(\xi)\} = \xi^{-2} v(\xi)^{d}$, we get
\[
\exp \bigl[ - C\xi^2 \psi\bigl\{v(\xi)\bigr\}
\bigr] = \exp \bigl\{- C v(\xi)^d \bigr\} = \exp\bigl\{- C
\xi^{-d/\beta} \log^{-d/2}(1/\xi) \bigr\}.
\]

Substituting this bound in (\ref{lbe56}), we finally obtain
%
%
\begin{equation}
\label{eqf1} P \bigl(\bigl\llVert W^A- \mu_0\bigr
\rrVert _2 \leq\xi \bigr) \leq \exp\bigl\{-C \xi^{-d/\beta}
\log^{-d/2}(1/\xi)\bigr\}.
\end{equation}
From \citet{van2009adaptive}, it also follows that
%
%
\begin{equation}
\label{eqf2} P\bigl(\bigl\llVert W^A - \mu_0\bigr
\rrVert _\infty< \varepsilon_n\bigr) \geq\exp \bigl\{- n
\varepsilon_n^2 \bigr\}
\end{equation}
for $\varepsilon_n = n^{-\alpha/(2\alpha+d)}\log^{t_1}n$. Recall
$\gamma_1 = \alpha/(2\alpha+ d)$. Using (\ref{eqf1}) and (\ref{eqf2}),
we obtain with $\xi$ replaced by $\xi_n = n^{-\gamma}$,
\[
\frac{P(\llVert  W^A - \mu_0\rrVert _2 < \xi_n)}{P(\llVert  W^A - \mu_0\rrVert _\infty< \varepsilon_n)} \leq\exp \bigl\{- \bigl(n^{\gamma d/\beta}
\log^{-d/2}n - n^{\gamma_1d / \alpha} \log ^{2t_1} n\bigr)\bigr\}.
\]
Since $\gamma/\beta> \gamma_1 /\alpha$ by assumption, $n^{\gamma
d/\beta} > 3n^{\gamma_1 d/\alpha} \log^{2t_1+ d/2} n$ for large\break
enough~$n$. Hence, for large enough $n$,
\[
\frac{P(\llVert  W^A - \mu_0\rrVert _2 < \xi_n)}{P(\llVert  W^A - \mu_0\rrVert _\infty< \varepsilon_n)} \leq\exp \bigl\{- 2 n^{\gamma_1d / \alpha} \log^{2t_1} n
\bigr\} = \exp\bigl\{- 2 n \varepsilon _n^2\bigr\},
\]
proving (\ref{eqlbmain}).

\section{Discussion}

We showed that a Gaussian process model with dimensional specific
scalings equipped with an appropriately chosen joint prior on the
scales can adapt to the true dimensionality or different smoothness
levels along different coordinates of the true function. In some
situations, it might be more reasonable to assume the true function to
be supported on a smaller dimensional linear subspace. In such cases, a
minor modification of our approach can achieve dimension adaptability
by incorporating an orthogonal projection of the covariate space as
$W^{\mathbf{a},Q}(t) = W(\mathbf{a}\cdot Qt)$ for a $d \times d$
orthogonal matrix $Q$. Such an approach is recently pursued by \citet{tokdar2011dimension} which assumes isotropy in the ambient dimensions
and uses the same prior as in Section~\ref{secadap}. One could easily
allow anisotropy in the rotated coordinate system using the unified
prior in Section~\ref{secconnec}. As a topic for future research, we
would like to explore consistent estimation of the dimension of the
true subspace and the subspace itself.

A salient feature of our prior (PD1)--(PD4) compared to \citet{zou2010nonparametric},
\citet{savitsky2011variable} is that the tail heaviness
of $A$ is related to the subset size of $S$. For larger subsets, the
tails of $A$ get lighter, resulting in down-weighted scalings for
larger subsets compared to smaller ones. It would be interesting to
explore the implied difference in practical performance from \citet{zou2010nonparametric}.

\begin{appendix}
\section{}\label{app0}
In this appendix, we provide examples of functions $\eta$ and $\zeta$
which satisfy condition (LBT) in Section~\ref{seclowerbounds}.

\begin{pf*}{Example of $\eta$}
We first provide an example of
$\eta$ with $u = 1/2$. Existence of \emph{bump functions} or
infinitely smooth\vspace*{1pt} compactly supported functions are well known; for
example, $\Psi(t) = \exp\{ -1/(1-t^2) \} \chi_{|t| < 1}$ is an
example of a $C^{\infty}$ function with support $[-1, 1]$ [Section~13
of \citet{tu2010introduction}]. From \citet{johnson2007saddle}, $|\hat
{\Psi}(\lambda)| \asymp\exp(-C\sqrt{|\lambda|})$ for large
$|\lambda|$. Define $\eta\dvtx  [0, 1] \to\mathbb{R}$ by $\eta(x) =
\Psi\{2(x-v)/(w-v) - 1\}$ if $x \in(v, w)$ and zero otherwise, that
is, shift and scale $\Psi$ to have support on $[v, w]$. A simple
calculation yields $\hat{\eta}(\lambda) = (2 \pi a)^{-1} e^{- i
\lambda/a} \hat{\Psi}(\lambda/a)$, where $a = 2/(w-v) - 1$ and $b =
-2v/(w-v)$. Hence, one also has $\llvert \hat{\eta}(\lambda
)\rrvert \asymp\exp(-C\sqrt{|\lambda|})$ for large $|\lambda
|$. Continue to denote by~$\eta$ the function on $[0, 1]^{d_2}$ given
by $\eta(x) = \prod_{j=1}^{d_2} \eta(x_j)$. Clearly, $\eta$ has
support $[v, w]^{d_2}$ and $|\hat{\eta}(\lambda)| \asymp\exp
(-C\sqrt{\llVert \lambda\rrVert })$, since\footnote{With
$a_j = | \lambda_j|^{1/2}$, $(\sum a_j)^4 \geq\sum a_j^4$ and using
Cauchy--Schwarz inequality twice, $(\sum a_j)^4 \leq d_2^{3} \sum a_j^4$.}
$d_2^{3/4} \llVert \lambda\rrVert ^{1/2} \geq\sum_{j=1}^{d_2} \sqrt{|\lambda_j|} \geq\llVert \lambda\rrVert
^{1/2}$. Transforming to polar coordinates,
\begin{eqnarray*}
\int_{\llVert \lambda\rrVert \in[K_2, 2K_2]} \bigl| \hat{\eta }(\lambda) \bigr|^2 \,d \lambda&=&
C' \int_{K_2}^{2 K_2} r^{d_2 -1}
e^{-2 C
\sqrt{r} } \,dr
\\
& =& C' \int_{\sqrt{K_2}}^{\sqrt{2K_2}}
z^{2 d_2-1} e^{-2 C z} \,dz \geq e^{-2 C \sqrt{K_2}}
\end{eqnarray*}
for large $K_2$.\noqed
\end{pf*}

\begin{pf*}{Example of $\zeta$}
As mentioned in the final paragraph of Section~\ref{seclowerbounds},
we present a concrete example of $\zeta$ in the case $d_1 =1$ with
$\alpha=1$ and $\beta= 3/2$. 
For ease of notation, we present the example on $[-2, 2]$ with support
$[-1, 1]$; linearly transforming to any compact interval which is a
subset of $[0, 1]$ does not affect the tail behavior of the Fourier transform.

Define $\zeta\dvtx  [-2, 2] \to\mathbb{R}$ as $\zeta(x) = (1- \llvert  x\rrvert ^2)\chi_{\llvert  x\rrvert \leq1}$.
Observe that $\zeta\in C^{1}[-2, 2]$, that is, $\zeta$ is Lipschitz
continuous, since $\zeta$ is absolutely continuous with an a.e.
bounded derivative. However, $\zeta\notin C^{1+s}[-2, 2]$ for any $s >
0$ as $\zeta$ is not differentiable at $\pm1$. By equation~9.1.20 in \citet{abramowitz1970handbook},
$\hat{\zeta}(\lambda) = \sqrt{2/\pi} \llvert \lambda\rrvert ^{-3/2} J_{3/2}(\llvert \lambda\rrvert )$, where
$J_{\nu}(x)$ is the Bessel function of the first kind of order $\nu$.
Further, combining equations 10.1.1 and 10.1.11 in \citet{abramowitz1970handbook},
$J_{3/2}(x) = \sqrt{2/(\pi x)} (\sin x/x - \cos x)$ for $x > 0$, so that
\[
\hat{\zeta}(\lambda) = \biggl(\frac{2}{\pi} \biggr) |\lambda|^{-3}
\bigl(\sin|\lambda| - |\lambda| \cos|\lambda|\bigr).
\]

Using $\int_{M}^{\infty} x^{-l} \cos^2(x) \,dx, \int_{M}^{\infty}
x^{-l} \sin^2(x) \,dx \asymp M^{-l +1}$ for any $l > 2$ and $M > 0$
large,\footnote{See Appendix~\ref{app3}.} we have for sufficiently large $K_1 > 0$,
\begin{eqnarray*}
&& \int_{\llvert \lambda\rrvert \geq K_1}\bigl|\hat{\zeta }(\lambda)\bigr|^2 \,d \lambda
\\
&&\qquad =
C_1 \biggl\{ \int_{K_1}^{\infty}
\lambda^{-6} \sin^2(\lambda)\,d \lambda+ \int
_{K_1}^{\infty} \lambda^{-4}
\cos^2(\lambda)\,d\lambda- \int_{K_1}^{\infty}
\lambda^{-5} \sin(2 \lambda)\,d \lambda \biggr\}
\\
&&\qquad \geq C_2 K_1^{-5} + C_3
K_1^{-3} - C_4 K_1^{-4}
\\
&&\qquad \geq C_5 K_1^{-3}
\end{eqnarray*}
for constants $C_i> 0, i=1, \ldots, 5$. Thus, (LBT) is satisfied with
$\beta= 3/2$.\noqed
\end{pf*}

\section{}\label{app1}
%
%
\begin{prop}\label{remexistr}
Let $0 < v < w <1$. There exists a function $r\dvtx  [0, 1]^d \to\mathbb
{R}$ such that $r(t) = 1$ for all $t \in[v, w]^d$, $r$ has support
within $[0, 1]^d$ and $\llvert \hat{r}(\lambda)\rrvert \leq
e^{-K \sqrt{\llVert \lambda\rrVert }}$ for large $\llVert \lambda\rrVert $.
\end{prop}
\begin{pf}
We first construct a $C^{\infty}$ function $\rho\dvtx  \mathbb{R} \to[0,
1]$ which is identically~1 on $[-a,a]$ and has support in $[-b, b]$ for
some $0 < a < b < 1$.
Recall the $C^{\infty}$ bump function $\Psi$ from Appendix~\ref{app0} and set $g(t) = (1/c)\Psi(t/c)$ with $c = (b-a)/2$. Clearly,
$g$ is an infinitely smooth function with support $[-c, c]$. Let $\xi\dvtx \mathbb{R} \to[0, 1]$ be the indicator function of the interval
$[-(b+a)/2, (b+a)/2]$. Define $\rho= g * \xi$. We claim that $\rho$
is a smooth function identically 1 on $[-a,a]$ and has support in $[-b,
b]$. Observe that if $\llvert  x\rrvert > b$, $\rho(y) = \int_{-\infty}^{\infty} g(y)\xi(x-y) \,dy =0$ as $\llvert  x-y\rrvert >(b+a)/2$ if $\llvert  y\rrvert < c = (b-a)/2$. Also,
if $\llvert  x\rrvert < a$, $\llvert  y\rrvert \leq
c$, $\llvert  x - y\rrvert \leq a+c =(a+b)/2$. Hence, $\rho
(x) = 1$ if $\llvert  x\rrvert < a$. It is also easy to show that
$\rho\dvtx  \mathbb{R} \to[0, 1]$ is a $C^{\infty}$ function.

Now, map the interval $[-a, a]$ linearly to $[v, w]$ through $x \mapsto
px + q$, with $p = (w-v)/(2a)$ and $q = (w+v)/2$. Also, let $w_1 = pb +
q$ and $v_1 = -pb + q$. By suitable choice of $a$ and $b$, we can
ensure that $0 < v_1 < w_1 < 1$. Then the function $\tilde{\rho}(t) =
\rho\{(t-q)/p\}$ is infinitely smooth, equals $1$ on $[v, w]$ and has
support in $[v_1, w_1]$. Defining $r(t) = \tilde{\rho}(t_1) \cdots
\tilde{\rho}(t_d)$, all assertions of Proposition~\ref{remexistr}
are satisfied barring the tail behavior of the Fourier transform which
we prove below. Proceeding as in Appendix~\ref{app0}, one has
\[
\bigl\llvert \hat{r}(\lambda)\bigr\rrvert = \Biggl| \prod
_{j=1}^d\hspace*{1.5pt} \hat{\hspace*{-1.5pt}\tilde{\rho}}(\lambda_j)
\Biggr| = C \prod_{j=1}^d \bigl| \hat{g}(p \lambda
_j)\bigr| \bigl| \hat{\xi}(p \lambda_j)\bigr|.
\]
Using $|\hat{\Psi}(\lambda)| \asymp\exp(-C\sqrt{|\lambda|})$ and
$|\hat{\xi}(\lambda)| \leq1$ for large $|\lambda_j|$, one has
$\llvert \hat{r}(\lambda)\rrvert \leq C\prod_{j=1}^d\exp
(-K_1\sqrt{|\lambda_j|})$ for constants $C, K_1>0 $ and for large
$\llvert \lambda_j\rrvert $. The proof is completed by
observing $\sum_{j=1}^d \sqrt{|\lambda_j|} \geq\llVert \lambda
\rrVert ^{1/2}$.\vspace*{-3pt}
\end{pf}

\section{}\label{app2}
This is a modified version of Lemma 16 in \citet{van2011information};
the proof is a simple extension, and hence omitted.
%
%
\begin{prop}\label{lemmodlem16}
For arbitrary functions $f,g\dvtx \mathbb{R}^d \to\mathbb{R}$, $I \subset
\{1, \ldots, d\}$, $J = \{1, \ldots, d\} \setminus I $; $\chi_{K_1,
K_2}$ and $\chi_{K_1,\cdot}$ the\vspace*{1pt} indicator functions of $\{\lambda
\in\mathbb{R}^d\dvtx  \llVert \lambda_{I}\rrVert > K_1, \llVert \lambda_{J}\rrVert \in[K_2, 2K_2] \}$ and $\{\lambda\in
\mathbb{R}^d\dvtx  \llVert \lambda_{I}\rrVert > K_1\}$,
respectively, and $0 <\break  A_1 < K_1$,
%
%
\begin{eqnarray}\label{eqL16}
&& \llVert f \chi_{K_1 - A_1, \cdot}\rrVert _{2, \mathbb{R}} \int
_{\llVert  t_I\rrVert \leq A_1} \bigl\llvert g(t)\bigr\rrvert \,dt
\nonumber\\[-9pt]\\[-9pt]
&&\qquad \geq \bigl\llVert(f* g) \chi_{K_1, K_2}\bigr\rrVert _{2,
\mathbb{R}}
-\llVert f\rrVert _{2, \mathbb
{R}} \int_{\llVert  t_I\rrVert > A_1} \bigl\llvert
g(t)\bigr\rrvert \,dt.\nonumber\vspace*{-3pt}
\end{eqnarray}
\end{prop}
%
\section{}\label{app3}
We show that for any $l > 2$ and $M > 0$ large, $\int_{M}^{\infty}
x^{-l} \sin^2(x) \,dx \asymp\break  M^{-l+1}$. The upper bound is immediate and
we focus on the lower bound. Without loss of generality, assume $M =
\pi m$ for some positive integer $m$, so that it is enough to consider
$\mathrm{I} = \int_{m}^{\infty} x^{-l} \sin^2(\pi x) \,dx$. Write
$\mathrm{I} =\break  \sum_{j=m}^{\infty} \int_{j}^{j+1} x^{-l} \sin^2(\pi
x) \,dx$. Noting that $\sin^2(\pi x)$ can be bounded below by $1/2$ on
$[j+1/4, j+3/4]$ for any $j \geq1$, we have $\mathrm{I} \gtrsim\sum_{j=m}^{\infty} b_j$, where $b_j = (j+1/4)^{-l+1} - (j+3/4)^{-l+1}$.
Noting that $(j+3/4)^{l-1} - (j+1/4)^{l-1} \geq(j+3/4)^{l-2}/2$, we
have $b_j \geq(j+3/4)^{-l}/2 \geq(2j)^{-l}/2$. Hence, $\sum_{j=m}^{\infty} b_j \gtrsim\sum_{j=m}^{\infty} \int_{j}^{j+1}
j^{-l} \,dx \geq\sum_{j=m}^{\infty} \int_{j}^{j+1} x^{-l} \,dx \geq
\int_{m}^{\infty} x^{-l} \,dx = m^{-l+1}/(l-1)$.

Along similar lines, we can show that $\int_{M}^{\infty} x^{-l} \cos
^2(x) \,dx \asymp M^{-l+1}$.\vspace*{-3pt}
\end{appendix}

\section*{Acknowledgments}
We sincerely thank the Associate Editor for suggesting the function $\zeta$
mentioned in Appendix~\ref{app0} which is used to demonstrate the
lower bound results. The authors also thank all the anonymous referees
for their detailed comments on previous versions of the paper which has
led to a correct proof of the lower bound result and also led to a
significantly improved presentation.


%

\printaddresses

\end{document}